\documentclass[12pt]{amsart}
\usepackage{graphicx}
\usepackage{times}
\usepackage{mathptm}

\begin{document}
\newtheorem{theorem}{Theorem}[section]
\newtheorem{prop}[theorem]{Proposition}
\newtheorem{lemma}[theorem]{Lemma}
\newtheorem{claim}[theorem]{Claim}
\newtheorem{cor}[theorem]{Corollary}
\newtheorem{defin}[theorem]{Definition}
\newtheorem{example}[theorem]{Example}
\newtheorem{xca}[theorem]{Exercise}
\newcommand{\map}{\mbox{$\rightarrow$}}
\newcommand{\aaa}{\mbox{$\alpha$}}
\newcommand{\Aaa}{\mbox{$\mathcal A$}}
\newcommand{\bbb}{\mbox{$\beta$}}
\newcommand{\ccc}{\mbox{$\mathcal C$}}
\newcommand{\ddd}{\mbox{$\delta$}} 
\newcommand{\Ddd}{\mbox{$\Delta$}}
\newcommand{\eee}{\mbox{$\epsilon$}}
\newcommand{\Fff}{\mbox{$\mathcal F$}}  
\newcommand{\Ggg}{\mbox{$\Gamma$}}
\newcommand{\ggg}{\mbox{$\gamma$}}
\newcommand{\kkk}{\mbox{$\kappa$}}
\newcommand{\lll}{\mbox{$\lambda$}}
\newcommand{\Lll}{\mbox{$\Lambda$}}
\newcommand{\rrr}{\mbox{$\rho$}} 
\newcommand{\sss}{\sigma} 
\newcommand{\Rbb}{\mbox{$\mathbb{R}$}} 
\newcommand{\Sss}{\mbox{$\mathcal S$}} 
\newcommand{\Ss}{\mbox{$\Sigma$}}
\newcommand{\Th}{\mbox{$\Theta$}} 
\newcommand{\ttt}{\mbox{$\tau$}} 
\newcommand{\bdd}{\mbox{$\partial$}}
\newcommand{\zzz}{\mbox{$\zeta$}}
\newcommand{\inter}{\mbox{${\rm int}$}}

\title[] {$3$-manifolds with planar presentations and the width of
satellite knots}

\author{Martin Scharlemann}
\address{\hskip-\parindent
        Mathematics Department\\
        University of California\\
        Santa Barbara, CA 93106\\
        USA}
\email{mgscharl@math.ucsb.edu}

\author{Jennifer Schultens}
\address{\hskip-\parindent
        Department of Mathematics and Computer Science\\
        Emory University\\
        Atlanta, Georgia 30322\\
        USA}
\email{jcs@mathcs.emory.edu}

\date{\today} \thanks{The authors thank RIMS Kyoto, where this work was
begun, Professor Tsuyoshi Kobayashi for inviting us to RIMS, Yo'av
Rieck for helpful conversations there, and NSF for partial support via
grants DMS 0203680 and DMS 0104039.  The second author also thanks the
MPIM-Bonn for support.}

\begin{abstract} We consider compact $3$-manifolds $M$ having a
submersion $h$ to $R$ in which each generic point inverse is a planar
surface.  The standard height function on a submanifold of $S^{3}$ is
a motivating example.  To $(M, h)$ we associate a connectivity graph
$\Ggg$.  For $M \subset S^{3}$, $\Ggg$ is a tree if and only if there
is a Fox reimbedding of $M$ which carries horizontal circles to a
complete collection of complementary meridian circles.  On the other
hand, if the connectivity graph of
$S^{3} - M$ is a tree, then there is a level-preserving reimbedding of
$M$ so that $S^{3} - M$ is a connected sum of handlebodies.

Corollary: 
\begin{itemize}
    
    \item The width of a satellite knot is no less than the width of its
    pattern knot and so
    
    \item $w(K_{1} \# K_{2}) \geq max(w(K_{1}),
    w(K_{2}))$
    
    \end{itemize}

\end{abstract}
\maketitle
  

The notion of thin position, introduced by D. Gabai \cite{G}, has been
employed with great success in many geometric constructions.  Yet the
underlying notion of the width of a knot remains shrouded in mystery. 
Little is known about the width of specific knots, or how knot width
behaves under connected sum.  By stacking a copy of $K_1$ in thin
position on top of a copy of $K_2$ in thin position, it is easily seen
that $w(K_1 \# K_2) \leq w(K_1) + w(K_2) -2$.  Here we establish a
lower bound for the width of a knot sum: the width is bounded below by the
maximum of the widths of its summands and therefore also by one half
the sum of the widths of its summands.

Knot width can be thought of as a kind of refinement of bridge number. 
Interest in how the width of a knot behaves under connected sum is
inspired, in part, by the fact that bridge number behaves very well. 
Indeed for bridge number, $b(K_1 \# K_2) = b(K_1) + b(K_2) - 1$, see
the paper \cite{S} by H. Schubert, or \cite{Sc3} for a much shorter
proof.  The shorter proof in \cite{Sc3} crystallized out of an
investigation into whether or not thin position arguments clarify the
behaviour of bridge number under connected sum.  The answer to that
question appears to be no: width seems to be a much more refined
invariant than can be useful for the recovery of Schubert's result. 
In particular, the argument in \cite{Sc3} fails in settings where the
swallow follow torus is too convoluted.  One suspects that
degeneration of width under connected sum of knots is possible, i.e.,
that there might be knots $K_1, K_2$, such that $w(K_1 \# K_2) <
w(K_1) + w(K_2) - 2$.  The situation may be analogous to that of
another knot invariant, tunnel number.  For small knots (knots whose
complements contain no essential closed surfaces), neither width nor
tunnel number degenerate under connected sum; i.e., for small knots,
width of knots satisfies $w(K_1 \# K_2) = w(K_1) + w(K_2) - 2$ and
tunnel number satisfies $t(K_1 \# K_2) = t(K_1) + t(K_2)$.  This is
proven in \cite{RS} and \cite{MS}, respectively.  On the other hand,
it is known that tunnel number can degenerate under connected sum, for
knots that are not small.  See for example \cite{Mo}.  Our results on
knot width are in a spirit similar to that of \cite{ScSc},
establishing an upper bound for such possible degeneration. 
Explicitly:

\bigskip

\noindent {\bf Corollary \ref{cor:sum}}
{\it For any two knots $K_1, K_2$, \[w(K_1 \# K_2) \geq max \{w(K_1),
w(K_2) \} \geq  \frac{1}{2}(w(K_1) + w(K_2)).\]}

\bigskip

We obtain Corollary \ref{cor:sum} by applying the following more
general result to the swallow-follow companion tori that are
associated to the connected sum of knots (see \cite[p.  10]{L}, or the
discussion in Section \ref{section:width}).

\bigskip

\noindent {\bf Corollary \ref{cor:satellite}} {\it Suppose $K'$ is a
satellite knot with pattern $K$.  Then $w(K') \geq w(K)$.}

\bigskip

Our approach to the latter result is to think of the companion solid
torus as a simple example of a handlebody in $S^{3}$.  We ask, in
general, how a handlebody $H$ in $S^{3}$ might be reimbedded so that
its complement is also a handlebody, hoping in particular to find a
reimbedding that preserves the natural projection to $R$ (called {\em
height}: $h: H \subset S^{3} \subset R^{4} \map R$).  There is a
theory of reimbeddings in $S^{3}$ going back to Fox \cite{Fo}, who
showed that any connected $M \subset S^{3}$ can be reimbedded so that
its complement is a union of handlebodies.  What is new here is the
concern about height $h: M \map R$.

In Section \ref{section:fox} we associate to an arbitrary compact $M
\subset S^{3}$, a certain graph $\Ggg$, and show that $\Ggg$ is a tree
if and only if there is a collection of horizontal (with respect to
height) circles in $\bdd M$ which constitute a complete collection
of meridian circles after a reimbedding whose complement is a
handlebody.  This discussion is in some sense only a digression; the
main argument begins with Section \ref{section:unknot}. 

Our goal is to reimbed a handlebody $H$ (preserving height) so that
the complement $M = S^{3} - H$ is also a handlebody.  What we in fact
study carefully is the complement $M$, hoping that by reconstructing
it appropriately, without changing $h$ on $M$, we can turn $M$ into a
handlebody.  One way to recognize that we are done is to observe that
if $H$ can be made to look like the neighborhood of a graph $\Lll$ and
$\Lll$ lies in $S^{2} \subset S^{3}$ then $S^{3} - H$ is indeed a
handlebody.  We call such a graph $\Lll$ unknotted.  In Section
\ref{section:unknot} we develop methods to construct and recognize
unknotted graphs.  In Section \ref{section:braid} we describe how, if
the graph $\Ggg$ associated to $M = S^{3} - H$ is a tree, we can
reconstruct $M$, without affecting height $h$, so that $M$ becomes the
complement of an unknotted graph, i.  e. a handlebody.  Such a
reimbedding of $H$ is called a Heegaard reimbedding.  In
Section \ref{section:heegaard} we observe that the only effect of this
reconstruction of $M$ on $H$ is to alter it by braid moves; the
corollaries on knot width then follow in Section \ref{section:width}.

\section{Mathematical preliminaries}

Throughout the paper, all manifolds will be orientable and, unless
otherwise stated, compact.  All embeddings will be locally flat. 
Since in dimension three there is little topological distinction between
smooth manifolds and PL manifolds, and it will be convenient to use
ideas and language from both smooth and PL topology, we will do so
without apology, leaving it to the reader to make the appropriate
translation if a specific structure (smooth or PL) is initially given
on the manifold.  

\begin{defin}  A {\em planar presentation} of a $3$-manifold 
$(M, \bdd M)$, $ \bdd M \neq \emptyset$ is a map $h: M \map 
R$ so that 
\begin{enumerate}
    \item $Dh: T_M \map T_{R}$ is always surjective
    
    \item $h|\bdd M$ is a general position Morse function and
    
    \item for $t$ any regular value of $h|\bdd M$, $h^{-1}(t)$ is a 
    planar surface, denoted $P^{t}$.
    
    \end{enumerate}
    
\end{defin}

The motivating source of examples is this: Consider $S^{3} \subset
R^{4}$ and let $p: R^{4} \map R$ be a standard projection, so
$p|S^{3}: S^{3}\map [-1, 1]$ has two critical points in $S^{3}$,
typically called the north and south poles.  Let $M \subset S^{3}$ be
a compact submanifold that does not contain either pole.  Then $h =
p|M$ is a planar presentation of $M$.  

\begin{figure}[tbh]
\centering
\includegraphics[width=.8\textwidth]{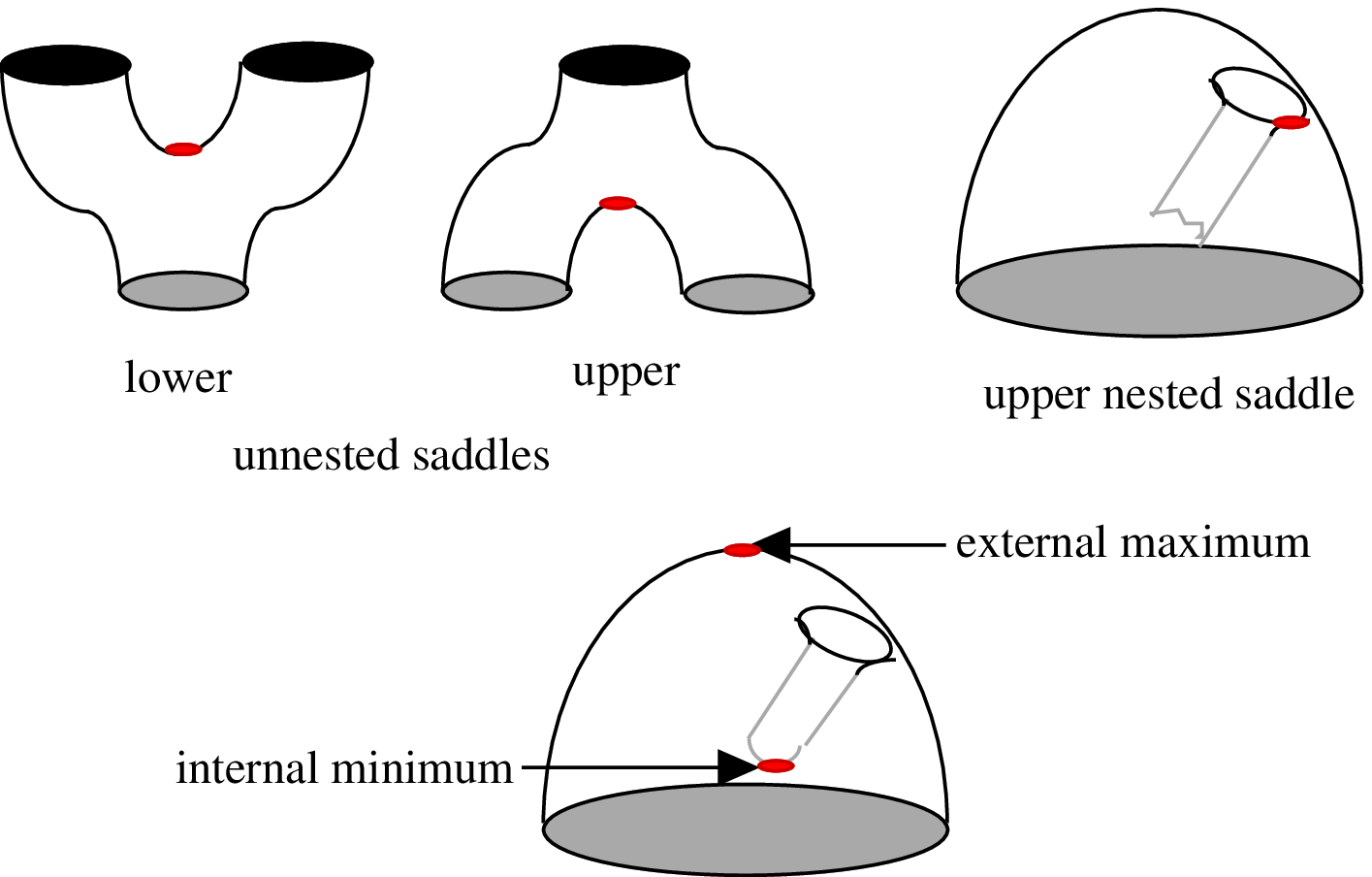}
\caption{} \label{fig:saddles}
\end{figure}

Consider an index one (i.  e. saddle) critical value $t_{\sigma}$
of $h|\bdd M$.  The corresponding critical point is called an {\em
upper saddle} (resp.  {\em lower saddle}) if $\bdd P^{t_{\sss}-
\epsilon}$ has one more (resp.  one less) circle component than $\bdd
P^{t_{\sss} + \epsilon}$.  If the number of components in $P^{t_{\sss}
+ \epsilon}$ and $P^{t_{\sss}- \epsilon}$ is the same, we say the
saddle is {\em nested}; otherwise the saddle is {\em unnested}.  Here
is an alternate description: an upper (resp.  lower) saddle is nested
if and only if the outward normal from $M$ points up (resp.  down) at
the saddle point.  (In particular, if $S \subset S^{3}$ is a surface,
then a saddle singularity of $S$ is a nested saddle for the component
of $S^{3} - S$ lying just above the singularity if and only if it is
unnested for the component of $S^{3} - S$ lying just below the
saddle.)  Similarly, a maximum (resp.  minimum) of $h$ on $\bdd M$ is
called an external maximum (resp.  minimum) of $h$ on $M$ if the
outward pointing normal from $M$ points up (resp.  down) at the
critical point.  Other maxima and minima on $\bdd H$ will be called
{\em internal} maxima and minima.  See Figure \ref{fig:saddles}.

\section{The connectivity graph and Fox reimbedding}
\label{section:fox}

Throughout this section, $(M, h)$ will be a planar presentation,
$s_{1} < s_{2} < \ldots < s_{n}$ will be the set of critical values at
which $h|\bdd M$ has an unnested saddle or an external maximum or 
minimum.  The points $x_{1}, \ldots, x_{n} \in \bdd M$ will be the
corresponding critical points.

\begin{lemma} \label{lemma:vertex} Suppose $M_{0}$ is a component of 
$M - \cup_{i=1}^{n} P^{s_{i}}$.  Then for any generic height 
$\overline{t}$, $P_{0} = M_{0} \cap P^{\overline{t}}$ is connected 
(possibly empty).
\end{lemma}

\begin{proof} Choose any two points in $P_{0}$.  Since $M_{0}$ is 
connected, there is an arc $\aaa \subset M_{0}$ that runs between 
them; a generic such arc will have its critical heights at different 
levels than $\bdd M$ does.  Since $\aaa \subset M_{0}$, $\aaa$ is 
disjoint from $\{ P^{s_{i}} \}$.  So for some $i$, the height of 
$\aaa$ lies between $s_{i}$ and $s_{i+1}$.  Let $t_{1}, \ldots, t_{m}$ 
be the critical values (if any) of $h|\bdd M_{0}$ between $s_{i}$ and 
$s_{i+1}$ and choose $\aaa$ to minimize the number of points in $T 
_{\aaa}= \{ t_{j} \in h(\aaa) \}$.  If $T_{\aaa}$ is empty, then 
$\aaa$ lies entirely in a region with no critical values, i.  e.  
$M_{0} \cap h^{-1}h(\aaa) \cong P_{0} \times I$.  Project $\aaa$ to 
$P_{0}$ and deduce that the ends of $P_{0}$ lie in the same component 
of $P_{0}$.

We now show that in fact $T_{\aaa}$ is always empty.  For suppose
$t_{j}$ is the greatest value (if any) of $T_{\aaa}$ that is greater
than $\overline{t}$ (or, symmetrically, the lowest value of $T_{\aaa}$
below $\overline{t}$).  The same argument as above shows that each
subarc of $\aaa$ that lies above $P^{t_{j}}$ can be projected to lie
in $P^{t_{j} + \epsilon}$ for any small $\epsilon$.  Since, by
assumption, passing through the critical level $t_{j}$ does not
connect or disconnect any component of $P^{t}$, in fact such a subarc
can then be pushed below $t_{j}$.  Once this is done for every subarc
of $\aaa$ above $t_{j}$, $T_{\aaa}$ is reduced by the removal of $t_{j}$, a
contradiction.

We have thereby shown that any two points in $P_{0}$ can be connected
by an arc in $P_{0}$, so $P_{0}$ is connected.
\end{proof}

\bigskip

\begin{defin} The {\em connectivity graph} $\Ggg$ of $(M, h)$ is the
graph whose vertices correspond to components of $M - \cup_{i=1}^{n}
P^{s_{i}}$ and whose edges correspond to components of $\cup_{i=1}^{n}
(P^{s_{i}} - x_{i})$.  An edge corresponding to a component $P_{0}$ of
$P^{s_{i}} - x_{i}$ has its ends at the vertices that correspond to
the components of $M - \cup_{i=1}^{n} P^{s_{i}}$ that lie just above
and below $P_{0}$.
\end{defin}

It is an old theorem of Fox \cite{Fo} that any compact connected
3-dimensional submanifold $M$ of $S^{3}$ can be reimbedded in $S^{3}$
so that the closure of $S^{3} - M$ is a union of handlebodies.  (This
theorem has recently been updated to include other non-Haken
$3$-manifolds \cite{Th}.)  As described above, let $p: S^{3} \map R$
be the standard height function and $M \subset S^{3}$ be a
$3$-manifold in general position with respect to $p$.  One can refine
Fox's question and ask if $M$ can be reimbedded in $S^{3}$ so that the
complement is a collection $H$ of handlebodies and, furthermore, each
horizontal circle in $\bdd M$ (that is each component of each generic
$\bdd P^{t}$) bounds a disk in $H$.  Put another way, the question is
whether a Fox reimbedding of $M$ can be found so that in the
complementary handlebodies a complete collection of meridian disks is
horizontal with respect to the original height function on $M$.

A first observation is that we may as well assume $M$ does not contain 
the poles.  For if $M$ contains the north pole, say, let $t$ be the 
highest critical value of $h = p|M$ on $\bdd M$, necessarily the image 
of a maximum on $\bdd M$.  Alter $M$ by simply removing the ball 
$h^{-1}(t - \eee, \infty)$.  The result does not contain the north pole and 
(after a tiny isotopy) is homeomorphic to $M$ via a homeomorphism that 
preserves the height function $h$ on $\bdd M$.  So, after this initial 
reimbedding, we may think of the pair $(M, h)$ as a planar 
presentation of $M$.

Then the answer is straightforward:

\begin{prop} \label{prop:fox}
    
    There is a collection of handlebodies $H$ so that $M \cup_{\bdd} H
    \equiv S^{3}$.  Moreover, there is a complete collection of
    meridian disks for $H$ whose boundaries are all horizontal (with
    respect to $h$) in $\bdd M$ if and only if the connectivity graph
    $\Ggg$ of $M$ is a tree.
    
    \end{prop}
    
\begin{proof}  The first claim is the central theorem of \cite{Fo}.  

The second claim follows from the central theorem of \cite{Sc1}.  This
says that a collection of $0$-framed curves $C \subset \bdd M$
contains a complete collection of meridians for some complementary
handlebody $H$ if and only if it has this property: Any properly
embedded surface $S$ in $M$ whose boundary is disjoint from $C$
separates $M$.

If $\Ggg$ is not a tree then some component $P_{0}$ of some $P^{t}$
is non-separating and clearly such a component can be made
disjoint from any finite collection of horizontal circles in $M$.  If
$M$ could be imbedded in $S^{3}$ so that the complement consisted of
handlebodies $H$ in which a complete collection of meridian boundaries
were horizontal with respect to $h$, then $P_{0}$ could be capped off
in $H$ by adding disks to $\bdd P_{0}$.  The result would be a
non-separating closed surface in $S^{3}$, and this of course is
impossible.

Conversely, suppose $C$ is a finite collection of horizontal circles
in $\bdd M$ chosen so large that any horizontal circle in $\bdd M$ is
parallel to an element of $C$ in $\bdd M$.  Suppose $S$ is a generic
non-separating properly embedded surface in $M$ with boundary disjoint
from $C$.  Let $\aaa$ be a generic simple closed curve in $M$ which
intersects $S$ in an odd number of points.  Choose such an $S$ to
minimize $|S \cap (\cup_{i=1}^{n} P^{s_{i}})|$, where, as above, $\{
s_{i} \}$ is the set of heights of the unnested saddles and of the
minima and maxima of $M$.  

The first observation is that in fact $S \cap (\cup_{i=1}^{n} 
P^{s_{i}}) = \emptyset$.  For otherwise, choose an innermost circle 
$c$ of intersection of $S$ with a component $P_{0}$ of $\cup_{i=1}^{n} 
P^{s_{i}}.$ Here ``innermost'' means that $c$ cuts off from $P_{0}$ a 
subplanar surface $P_{-}$ whose boundary consists of $c$ and a 
collection of boundary circles of $P_{0}$.  Then replacing a vertical 
collar of $c$ in $S$ with two parallel horizontal copies of $P_{-}$ 
gives a surface which has fewer components of intersection with 
$\cup_{i=1}^{n} P^{s_{i}}$ but which still contains a non-separating 
component, since the number of intersections with $\aaa$ is increased 
by $2 \cdot |\aaa \cap P_{-}|$ and so remains odd.  Since the boundary 
of $P_{-}$ is horizontal, generically it is disjoint from $C$.

So $S$ lies in a component of $M - \cup_{i=1}^{n} P^{s_{i}}$ whose
closure we denote $M_{0}$.  Let $h(M_{0}) = [s_{i}, s_{i+1}]$, so
$M_{0}$ lies in a slice of $S^{3}$ homeomorphic to $S^{2} \times
[s_{i}, s_{i+1}]$.  So as not to be distracted by other parts of $M$,
let $Q$ be a $2$-sphere and momentarily think of $M_{0}$ as lying in
$Q \times [s_{i}, s_{i+1}]$.  Since every horizontal cross-section of
$M_{0}$ is connected, at any generic height a cross-section of $Q -
M_{0}$ is a collection of disks.  In particular, the boundary
components of $S$ can be capped off in $Q \times [s_{i}, s_{i+1}]$ to
give a closed surface $S_{+} \subset Q \times [s_{i}, s_{i+1}]$.  

Now consider how the arcs $\aaa \cap M_{0}$ lie in $Q \times [s_{i}, 
s_{i+1}]$.  Any arc with both ends in $Q \times (s_{i})$ or both 
ends in $Q \times (s_{i+1})$ can be entirely homotoped in $Q 
\times [s_{i}, s_{i+1}]$ into that end and so be made disjoint from 
$S_{+}$.  It follows that such an arc intersects $S$ an even number of 
times.  Since $\aaa$ intersects $S$ an odd number of times, it follows 
that there are an odd number of arcs of $\aaa \cap M_{0}$ that run 
from the top of $M_{0}$ to the bottom.  Then, returning again to $M 
\subset S^{3}$ there must be an odd number of arcs of $\aaa - M_{0}$ 
that run from the top of $M_{0}$ to the bottom of $M_{0}$ in $M - 
M_{0}$.  In particular, there is at least one such arc, so one can 
construct a closed curve in $M$ that intersects the bottom of $M_{0}$ 
in a single point $p$.  Hence removing the edge in $\Ggg$ 
corresponding to the component of $P^{s_{i}} - x_{i}$ in which $p$ 
lies does not disconnect $\Ggg$.  Since we can remove an edge and not 
disconnect $\Ggg$, $\Ggg$ is not a tree.
\end{proof}

\section{Unknotted graph complements} \label{section:unknot}

\begin{defin} For $N$ a compact $3$-manifold and $\Lll$ a finite
graph, a {\em proper embedding} $\Lll \subset N$ is an embedding so
that $\bdd N \cap \Lll$ consists of a collection of valence one
vertices of $\Lll$.  These vertices are denoted $\bdd \Lll$.  The
other vertices, some of which may also have valence one, are called
{\em interior vertices}.

In case $N = B^{3}, S^{3}$ or $S^{2} \times I$, the pair $(N - 
\eta(\Lll), \bdd N - \eta(\Lll))$ will be denoted $(N_{\Lll}, 
P_{\Lll})$ and will be called a {\em graph complement with planar part 
$P_{\Lll}$}.  Graphs $\Lll$ and $\Lll'$ are equivalent if there is a 
homeomorphism $(N, \eta(\Lll)) \cong (N, \eta(\Lll'))$.  In 
particular, if $\Lll'$ is any graph obtained from $\Lll$ by sliding 
and isotoping edges rel $\bdd \Lll$ then $\Lll$ and $\Lll'$ are 
equivalent graphs.  

Two graph complements $(N_{\Lll}, P_{\Lll})$ and $(N_{\Lll'},
P_{\Lll'})$ will be called equivalent if they are pairwise
homeomorphic.  In particular, equivalent graphs have equivalent graph
complements.
\end{defin}

\begin{lemma} \label{lemma:graphcompl} Suppose $(M, h)$ is a planar
presentation of a compact manifold $M$, $J \subset R$ is a proper generic
interval, and suppose that in a component $M^{J}$ of $h^{-1}(J)$ all
saddles are nested.  Then $M^{J}$ is homeomorphic to a graph
complement with planar part $h^{-1}(\bdd J) \cap M^{J}.$
\end{lemma}

\begin{proof} Since $M^{J}$ is connected and contains no unnested
saddles, every generic horizontal cross-section is a connected planar
surface, by Lemma \ref{lemma:vertex}.  Since $J$ is proper, $M^{J}$ is
compact, so we may as well assume $J$ is compact (say $J = [0, 1])$,
though we do not know that $h(M^{J}) = J$.  We will describe the graph
$\Lll$ for which $M^{J}$ is the complement; the details of the
homeomorphism then follow from standard Morse theory.  See Figure
\ref{fig:graph}.

\begin{figure}[tbh]
\centering
\includegraphics[width=.8\textwidth]{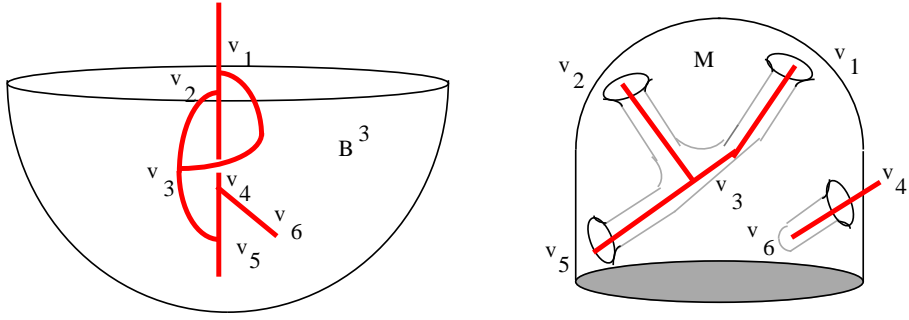}
\caption{$M$ as graph complement (in $B^{3}$)} \label{fig:graph}
\end{figure}

Suppose first that $J = h(M^{J})$, so each component of $h^{-1}(\bdd
J) \cap M^{J}$ is a non-empty connected planar surface.  We will
describe $\Lll \subset S^{2} \times I$.  Each (circle) boundary
component of $h^{-1}(\{ 1 \}) \cap M^{J}$ can be capped off with a
disk to give a two sphere; dually, $h^{-1}(\{ 1 \}) \cap M^{J}$ can be
thought of as obtained from $S^{2} \times \{1 \}$ by removing some
vertices.  These will be vertices in $\bdd \Lll$.  As $t$ descends
from $1$ through generic values of $t$, each boundary component of
$P^{t}$ can be capped off by a disk to give a sphere $S^{t}$.  This
gives an embedding $P^{t} \subset S^{t}$; dually $P^{t}$ can be
obtained from the sphere $S^{t}$ by removing a neighborhood of the
center of each disk.  As $t$ varies, these points form vertical edges
in $\Lll$ incident to those vertices of $\bdd \Lll$ that lie at height
$1$.

Now consider what happens as $t$ descends through a critical point of
$h|\bdd M$.  Each such critical point corresponds to an interior
vertex of $\Lll$.  In particular, edges descend from those vertices
that correspond to internal maxima (and descend {\em to} the vertices
that correspond to internal minima).  At lower saddles two edges
descend into the corresponding vertex and one edge descends from it
whereas at upper saddles one edge descends into the corresponding
vertex and two edges descend from it.  There are no external maxima or
minima, for these would necessarily start (or end) a different planar
surface, which could never be connected to $M^{J}$ since all saddles
are nested.

The argument is little changed if $h(M_{J}) \neq J$.  Say $1 \notin 
h(M^{J})$; then there is an external maximum on $M^{J}$ at height 
$t_{max} \in J$ and at a generic height just below it the ball $M^{J} 
\cap [t_{max}, 1]$ can be thought of as the complement of a radial arc 
in $B^{3}$, and so as a graph complement in $B^{3}$.  The rest of the 
construction proceeds as above, though now viewed as a construction in 
a collar $\bdd B^{3} \times I$.  Ultimately $M^{J}$ is thereby 
described as a graph complement in $B^{3}$ (or in $S^{3}$ if also $0 
\notin h(M^{J})$).
\end{proof}

The case when the graph $\Lll$ is planar will be particularly
important.  Let $S^{1} \times I$ denote the standard vertical cylinder
in $S^{2} \times I$.

\begin{defin} A properly imbedded graph $\Lll$ in $N = S^{3}$ (resp.  
$B^{3}$ or $S^{2} \times I$) is {\em unknotted} if it lies in $S^{2} 
\subset S^{3}$ (resp.  $B^{2} \subset B^{3}$ or $S^{1} \times I 
\subset S^{2} \times I$).  The pair $(M_{\Lll}, P_{\Lll})$ is then called 
an {\em unknotted graph complement} with {\em planar part $P_{\Lll}$}.  
(Note that the number of components of $P_{\Lll}$ determines whether 
the ambient manifold is $S^{3}, B^{3}$ or $S^{2} \times I$.)

More generally, any graph which is equivalent to an unknotted graph
will be called an unknotted graph.
\end{defin}

Unknotted graphs are in some sense unique:

\begin{prop} \label{prop:unique} Suppose $\Lll$ and $\Lll'$ are 
unknotted graphs in $N = S^{3}, B^3$, or $S^{2} \times I$.  Suppose 
that $\bdd \Lll = \bdd \Lll' \subset \bdd N$.  

Then there is a homeomorphism of pairs $(N, \eta({\Lll})) \cong (N, 
\eta({\Lll'})) $ which is the identity on $\eta(\bdd \Lll) = \bdd 
\eta(\bdd \Lll')$ if and only if the partitions of $\bdd \Lll = \bdd 
\Lll'$ defined by the components of $\Lll$ and $\Lll'$ are the same, 
and each component of $\Lll$ has the same Euler characteristic as the 
corresponding component of $\Lll'$.
\end{prop}

\begin{figure}[tbh]
\centering
\includegraphics[width=.8\textwidth]{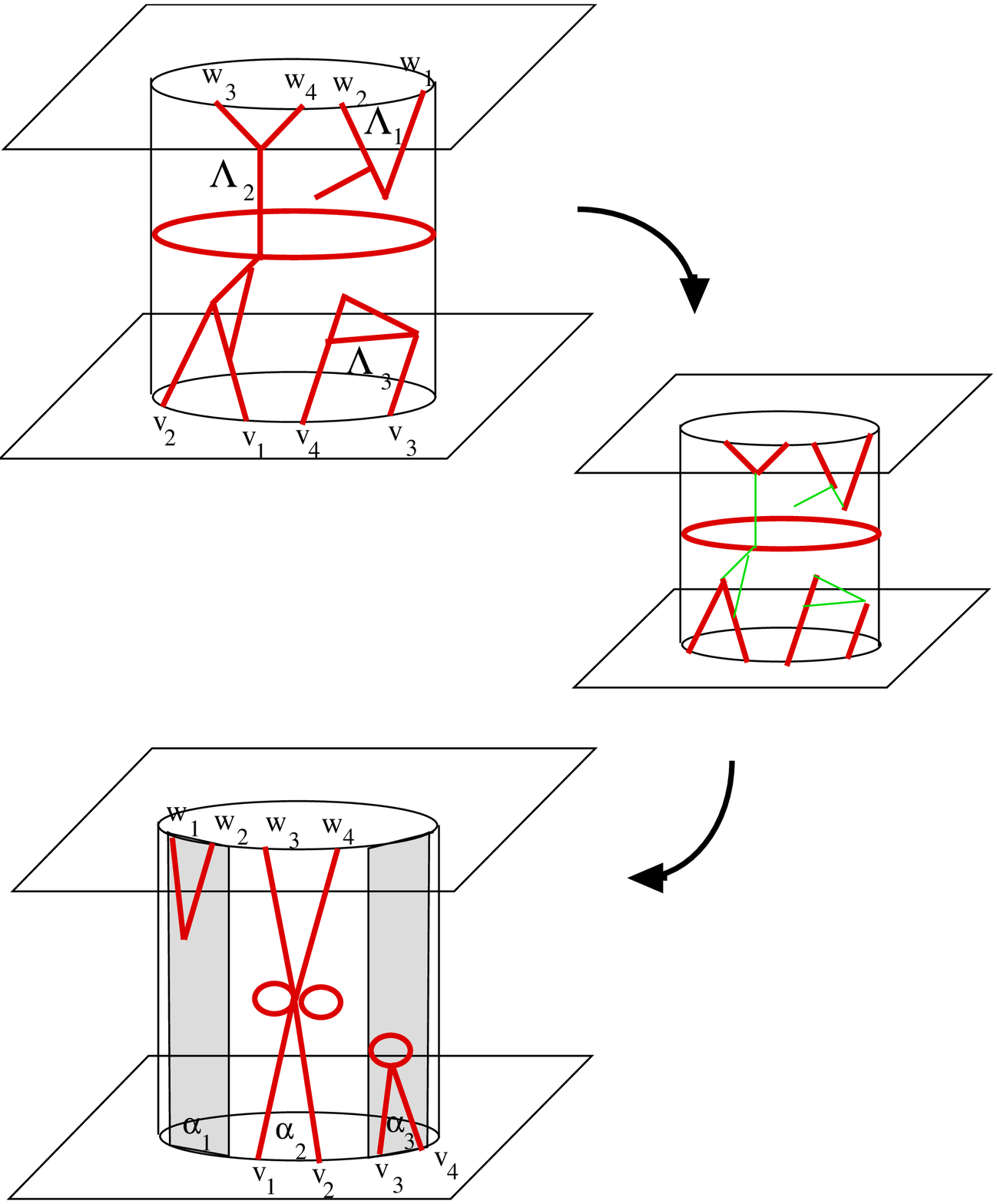}
\caption{Reimbedding $\Lll$ canonically} \label{fig:unique}
\end{figure}

\begin{proof} The existence of such a homeomorphism clearly implies 
that the partitions and the corresponding Euler characteristics are 
the same.  The difficulty is in proving the other direction.  We 
consider the case $N =S^{2} \times I$, for it is representative (and 
in fact the most difficult).

It will be convenient to number the $p$ components of $\Lll$ (and the
corresponding components of $\Lll'$) in some order $\Lll_{i}, i =
1,\ldots, p$, and then order the points $\bdd \Lll \cap (S^{2} \times
\{ 1 \}) = \{ w_{j} \}$ and $\bdd \Lll \cap (S^{2} \times \{ 0 \}) =
\{ v_{k} \}$ in some subordinate order, i.  e. so that in the ordering
all the boundary points of any earlier component of $\Lll$ come before
all the boundary points of any later component.

Since $\Lll$ is an unknotted graph we can assume (up to homeomorphism
of the pair $(N, \eta({\Lll}))$ rel $\eta(\bdd(\Lll))$) that $\Lll
\subset S^{1} \times I \subset S^{2} \times I$.  In a small
neighborhood of $\Lll$ collapse a forest that is maximal in $\Lll$
among those not incident to $\bdd \Lll$.  Then each component is the
cone on its boundary vertices, wedged with some circles.  Each circle
(even those that are essential in the vertical cylinder $S^{1} \times
I$) can be moved (rel the cone point) in $S^{2} \times I$ until it
bounds a tiny disk in $S^{1} \times I$ whose interior is disjoint from
$\Lll$.  For the purposes of the following argument, these tiny
circles can be ignored, since the assumption on Euler characteristic
means there will be as many tiny circles on a component of $\Lll'$ as
there are on the corresponding component of $\Lll$ (namely $1 -
\chi$).  Hence, with no loss of generality, we may assume $\Lll$ (and
$\Lll'$) consist entirely of collections of cones on (corresponding)
subsets of vertices.  See Figure \ref{fig:unique}.

Since $\Lll$ contains no circles there is a spanning arc of the
cylinder $S^{1} \times I$ that is disjoint from $\Lll$.  After an
isotopy in $S^{1} \times I$, we may as well assume the arc is vertical
and then break up a neighborhood of this vertical arc into a sequence
of $p$ vertical strips $\aaa_{i} \times I \subset S^{1} \times I, i =
1,\ldots,p$, where each $\aaa_{i} \cap \aaa_{i+1}, i = 1,\ldots, p-1$
is a single end point of each.  

Now push the first component $\Lll_{1}$ of $\Lll$ into a vertical 
cylinder $C$ parallel to $S^{1} \times I$ and, exploiting the fact 
that $\Lll_{1}$ is just a cone on its end points, do this so that the 
vertices in $\bdd \Lll_{1}$ appear in their correct order in a 
vertical strip in $C$.  Now move this vertical strip (and so 
$\Lll_{1}$) to $S^{1} \times I$ by moving the strip to $\aaa_{1} 
\times I$.  Similarly place the second component $\Lll_{2}$ in the 
second strip $\aaa_{2} \times I$ and continue through all of $\Lll$.  
Call the resulting graph $\Lll^{c(anonical)} \subset S^{1} \times I$ 
and observe that the process we have described gives a homeomorphism 
of pairs $g: (N, \eta({\Lll})) \map (N, \eta({\Lll^{c}})) $.  
Finally, observe that the process is so canonical that if we had done 
the same process on $\Lll'$ we would have obtained a homeomorphism of 
pairs $g': (N, \eta({\Lll'})) \map (N, \eta({\Lll^{c}}))$ 
that preserves the orderings.  In particular $g$ and $g'$ could be 
taken to be the same on $\eta (\bdd \Lll)) = \eta (\bdd 
\Lll'))$.  Then $g^{-1}g$ is the required homeomorphism of pairs.
\end{proof} 

\bigskip

If $\Lll$ is unknotted, then $M_{\Lll}$ has a particularly simple 
structure:

\begin{lemma} An unknotted graph complement is a connected sum of handlebodies.
\end{lemma}

Note: Here we regard a $3$-ball as a handlebody of genus $0$.

\begin{proof} The case in which the ambient manifold is $B^{3}$ is
representative.  We have $\Lll \subset B^{2} \subset B^{3}$.  In a
small neighborhood of $\Lll$ collapse a forest that is maximal among
all forests not incident to $\bdd \Lll$.  Then each component is the
cone on its boundary vertices, wedged with some circles.  Each circle
can be pushed out of $B^{2}$ (rel its wedge point) and made to bound a
tiny disk.  Removing such a circle from $\Lll$ has the effect in the
graph complement of removing a $1$-handle, dual to the tiny disk.  In
particular, with no loss, we can assume that no such circles arise and
so each component of $\Lll$ is a cone on its boundary vertices.

The proof is then by induction on $|\bdd \Lll|$.  If $\bdd \Lll = 
\emptyset$ then $\Lll$ is a collection of isolated vertices, so its 
complement is a connected sum of balls.  If any component of $\Lll$ 
has a single boundary vertex, then that component is just an arc with 
one end on $\bdd B^{3}$; removing it from $\Lll$ has no effect on the 
complement in $B^{3}$.  So without loss assume each component of 
$\Lll$ is the cone on two or more points in $\bdd B^{2}$.  A path in 
$\Lll$ between two such points divides the disk $B^{2}$ into two 
disks.  An outermost such path will cut off a disk $D$ from $B^{2}$ 
whose interior is disjoint from $\Lll$.  The disk $D$ can be used to 
$\bdd$-reduce $M_{\Lll}$ and the effect on $M_{\Lll}$ is the same as 
if we had removed one of the edges of $\Lll$ incident to $D$.  The 
proof then follows by induction.
\end{proof}

\bigskip

We now describe a few situations that guarantee that a graph is
unknotted in $S^{2} \times I$.  We will be taking the standard height
function on $S^{2} \times I$, namely projection to $I$.  A vertex $v$
in a properly embedded graph $\Lll \subset S \times I$ is a {\em
$Y$-vertex} if two or more edges are incident to $v$ from above and a
{\em $\lll$-vertex} if two or more edges are incident to $v$ from
below.  (A vertex may be both a $\lll$-vertex and a $Y$-vertex, or
neither.)


\begin{example} \label{exam:smallexam} Suppose $\Lll \subset S^{2}
\times I$ is a properly embedded graph so that

\begin{enumerate} 
  
    \item the edges in $\Lll$ are all monotonic with
respect to the projection $S^{2} \times I \map I$. 

\item there are no $Y$-vertices.

\end{enumerate}

Then $\Lll$ is an unknotted graph.
\end{example}

\begin{proof} We first simplify $\Lll$ up to graph equivalence.  By a
small edge-slide arrange that each vertex is incident to at most two
edges below; any vertex that is incident to a single edge above and a
single edge below can be ignored.  If an interior vertex is incident
to two edges below, and none above, then add a small vertical edge
above.  After these initial maneuvers, each interior vertex of $\Lll$
has valence zero, one or three; in the last case, the vertex is a
$\lll$-vertex.

Pick a circle $C$ in $S^{2} \times \{ 1 \}$ that contains all the 
vertices of $\bdd{\Lll}$ that lie in $S^{2} \times \{ 1 \}$.  As $t 
\in [0, 1]$ descends, the monotonicity of edges means that, until 
another vertex of $\Lll$ is encountered, the cross-section $\Lll \cap 
(S^{2} \times \{ t \})$ is a collection of points moving by isotopy in 
$S^{2}$.  Extend the isotopy to all of $S^{2}$ to get a continuously 
varying circle $C_{t} \subset S^{2} \times \{ t \}$ that contains all 
of $\Lll \cap (S^{2} \times \{ t \})$.  When a valence one vertex (or 
an isolated vertex) is encountered, it can be easily added to 
or deleted from $C_{t}$, as appropriate, depending on whether the edge 
incident to the vertex is incident from below or from above.

So we only need to worry about $\lll$-vertices.  As $t$ passes through
the level of such a vertex (which we have arranged to lie in $C_{t}$),
a single point in $C_{t}$ simply splits in two and we may incorporate
the arc between the two points as part of $C$.  Continue the process
down to $t = 0$.  Now, in a standard argument, the continuously
varying family of circles $C_{t}$ bounds a continuously varying family
of disks in $S^{2}$ and so there is a height-preserving isotopy of
$C_{t}$ to the standard $S^{1} \times I$.
\end{proof}

The fact that, in the proof, the original circle $C$ was ours to 
choose immediately leads to these additional examples:

\begin{example} \label{exam:concave} Suppose $\Lll \subset S^{2}
\times I$ is a properly embedded graph and there is a generic height 
$t \in I$ so that
\begin{enumerate}     
    \item the edges in $\Lll$ are all monotonic with
respect to the projection $S^{2} \times I \map I$.  

\item There are no $\lll$-vertices above $t$

\item There are no $Y$-vertices below $t$

\end{enumerate}

Then $\Lll$ is an unknotted graph.
\end{example}

\begin{figure}[tbh] 
\centering
\includegraphics[width=.5\textwidth]{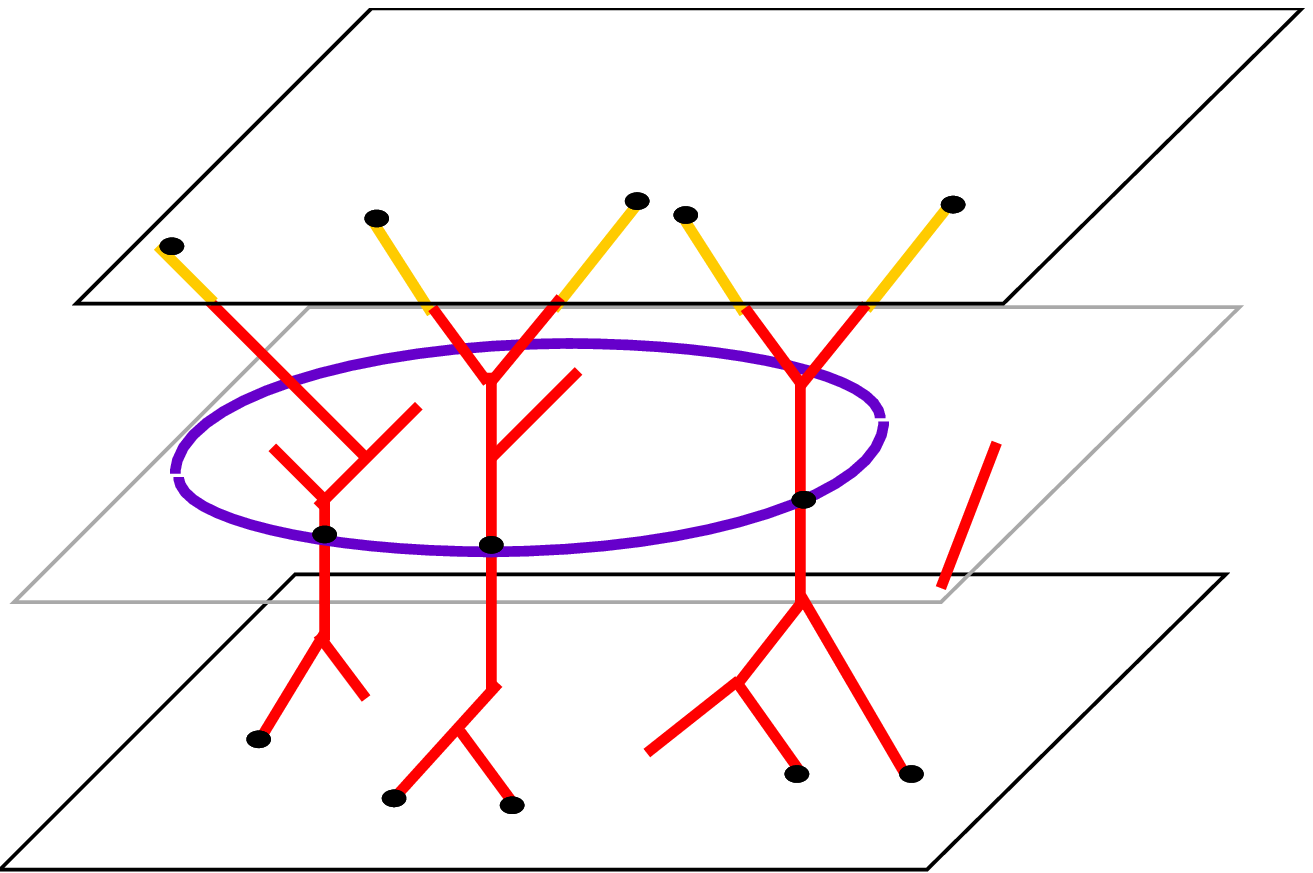}
\caption{} \label{fig:concave}
\end{figure}

\begin{proof} Apply the argument of Example \ref{exam:smallexam} 
separately to $S^{2} \times [0, t]$ and (upside down) to $S^{2} 
\times [t, 1]$, starting with a circle in $S^{2} \times \{ t \}$ that 
contains all points in $\Lll \cap (S^{2} \times \{ t \}).$  See 
Figure \ref{fig:concave}.
\end{proof}

More generally

\begin{example} \label{exam:pileon} Suppose $\Lll \subset S^{2}
\times I$ is a properly embedded graph so that
\begin{enumerate} 
    \item there is a generic level sphere $S^{2} \times \{ t \}$ for which 
    $\Lll$ intersects $S^{2} \times [t, 1]$ in an unknotted graph
    
    \item the edges in $\Lll \cap [0, t]$ are all monotonic with
respect to the projection to $[0, t]$   

    \item There are no $Y$-vertices in $\Lll \cap 
    [0, t]$.

\end{enumerate}

Then $\Lll$ is an unknotted graph.
\end{example}

\begin{proof} Apply the argument of Example \ref{exam:smallexam} 
separately to $S^{2} \times [0, t]$ starting with the circle in $S^{2} 
\times \{ t \}$ which is the base of the vertical cylinder in $S^{2} 
\times [t, 1]$ on which $\Lll \cap (S^{2} \times [t, 1])$ lies. \end{proof}

The next two examples simply reinterpret earlier examples in light of 
Lemma \ref{lemma:graphcompl}.

\begin{example} \label{exam:medexam} Suppose $(M, h)$ is a planar 
presentation of a manifold $M$ and for an interval $J$, $M^{J}$ is a 
component of $h^{-1}(J)$.  Suppose all saddles in $M^{J}$ are nested 
and that all lower saddles occur at higher levels than all the upper 
saddles do.  Then the pair $(M^{J}, P^{\bdd J})$ is an unknotted graph 
complement.
\end{example}

Suppose $(M, h)$ is a planar presentation of a 
manifold $M$ and $M^{[a,c]}$ is a component of $h^{-1}([a, c])$.  We 
will use the following notation: for $J$ a subinterval of $[a, c]$ let 
$M^{J} = M^{[a,c]} \cap h^{-1}(J)$ and for $t \in [a, c]$ let $Q^{t} = 
M^{[a,c]} \cap P^{t}$.

\begin{example} \label{exam:pileon2} Suppose $(M, h)$ is a planar 
presentation of a manifold $M$ and $M^{[a,c]}$ is a 
component of $h^{-1}([a, c])$.  
 
Suppose that for some $b \in [a, c]$

\begin{enumerate}

     \item $Q^{b}$ is connected
    
    \item the pair $(M^{[b,c]}, Q^{b} \cup Q^{c})$ is an unknotted 
    graph complement.
    
    \item all saddles in  $M^{[a, b]}$ are nested upper saddles. 

\end{enumerate}

Then $M^{[a, c]}$ is an unknotted graph complement.
\end{example}

It would be useful to know that if $\Lll_{1} \subset B_{1}$, $\Lll_{2}
\subset B_{2}$ are unknotted graphs in $3$-balls $B_{i}$, and we are
given some identification of $\bdd \Lll_{1}$ with $\bdd \Lll_{2}$,
that we could find some way to attach $\bdd B_{1}$ to $\bdd B_{2}$
consistent with that identification so that the resulting graph is
unknotted.  Ultimately we will succeed (see Lemma \ref{lemma:sum}) but
first we observe that the most obvious way to try to prove this fact is
doomed to fail.  Specifically, it may be impossible to match up the
boundary of a disk in $B_{1}$ containing $\Lll_{1}$ to the boundary of
a disk in $B_{2}$ containing $\Lll_{2}$ in a manner that preserves the
identification of $\bdd \Lll_{1}$ with $\bdd \Lll_{2}$.

To see that this is impossible, take the following simple example: let
each $\Lll_{i}$ be three copies of a cone on three points, so that
$\bdd \Lll_{i}$ is nine points.  Identify $\bdd \Lll_{1}$ with $\bdd
\Lll_{2}$ so that the resulting graph is the complete bipartite graph
$K_{3,3}$.  If one could identify the boundary of a disk containing
$\Lll_{1}$ with the boundary of a disk containing $\Lll_{2}$ in a way
consistent with the identification of $\bdd \Lll_{1}$ with $\bdd
\Lll_{2}$, we would have found an embedding of $K_{3,3}$ into the
$2$-sphere, which is famously impossible.

Yet there is a way to attach $\bdd B_{1}$ to $\bdd B_{2}$ so that
$\Lll_{1} \cup_{\bdd} \Lll_{2}$ is an unknotted embedding of $K_{3,3}$
in $S^{3}$; the argument above merely shows that, in order to
demonstrate that such an embedding is unknotted, edges will need to be
slid over edges, inevitably across the sphere $\bdd B_{i}$.  In other
words, the demonstration that there is an unknotted embedding of
$K_{3,3}$ is inevitably a bit harder than one might at first expect.  

It will be extremely useful to demonstrate that {\em any} bipartite
graph has an unknotted embedding in $S^{3}$, via a construction much
as above.  That is the goal of the following lemma.  Recall that a
bipartite graph with vertex sets $A$ and $B$ is a graph so that each
edge has one end among the vertices of $A$ and the other end among the
vertices of $B$.  We will show that any bipartite graph can be
imbedded in a very controlled way into a cube so that the embedded
graph is unknotted: that is, {\em after some edge slides} the graph
can be made to lie in a plane.  Some details of its structure will be
crucial in the discussion of braid equivalence in Section 
\ref{section:braid}.

\begin{lemma} \label{lemma:bipartite} Let $\Lll$ be a finite bipartite 
graph, with vertex sets $A$ and $B$.  Then there is an embedding 
of $\Lll$ in the cube $I \times [-1, 1] \times I$ so that:

\begin{enumerate}

\item $A = \{(i/|A|, -1, 0), i = 0, \ldots |A| - 1 \}$

\item $B = \{(j/|B|, 1, 0), j = 0, \ldots  |B| - 1 \}$

\item Each edge in $\Lll$ is monotonic with respect to the 
$y$-coordinate.  That is, each edge projects to $[-1, 1]$ 
with no critical points.

\item The edges may be isotoped and slid over each other (perhaps 
destroying the bipartite structure) in the cube, so that afterwards 
the resulting graph lies entirely in the face $I \times [-1, 1] \times 
\{ 0 \}$.

\end{enumerate}

Moreover, given a specific edge $e$ in $\Lll$, such an embedding can 
be found so that $e =\{0 \} \times [-1, 1] \times \{0 \}$ and $e$ 
never moves during the isotopy.

\end{lemma}

Note that the last numbered condition implies that $\Lll$ is an 
unknotted graph.  (Technically, $\Lll$ is unknotted only in a larger 
cube, for the given cube contains $e$ in a face and so does not 
contain $\Lll$ as a proper subgraph.)  

\begin{figure}[tbh] 
\centering
\includegraphics[width=1.0\textwidth]{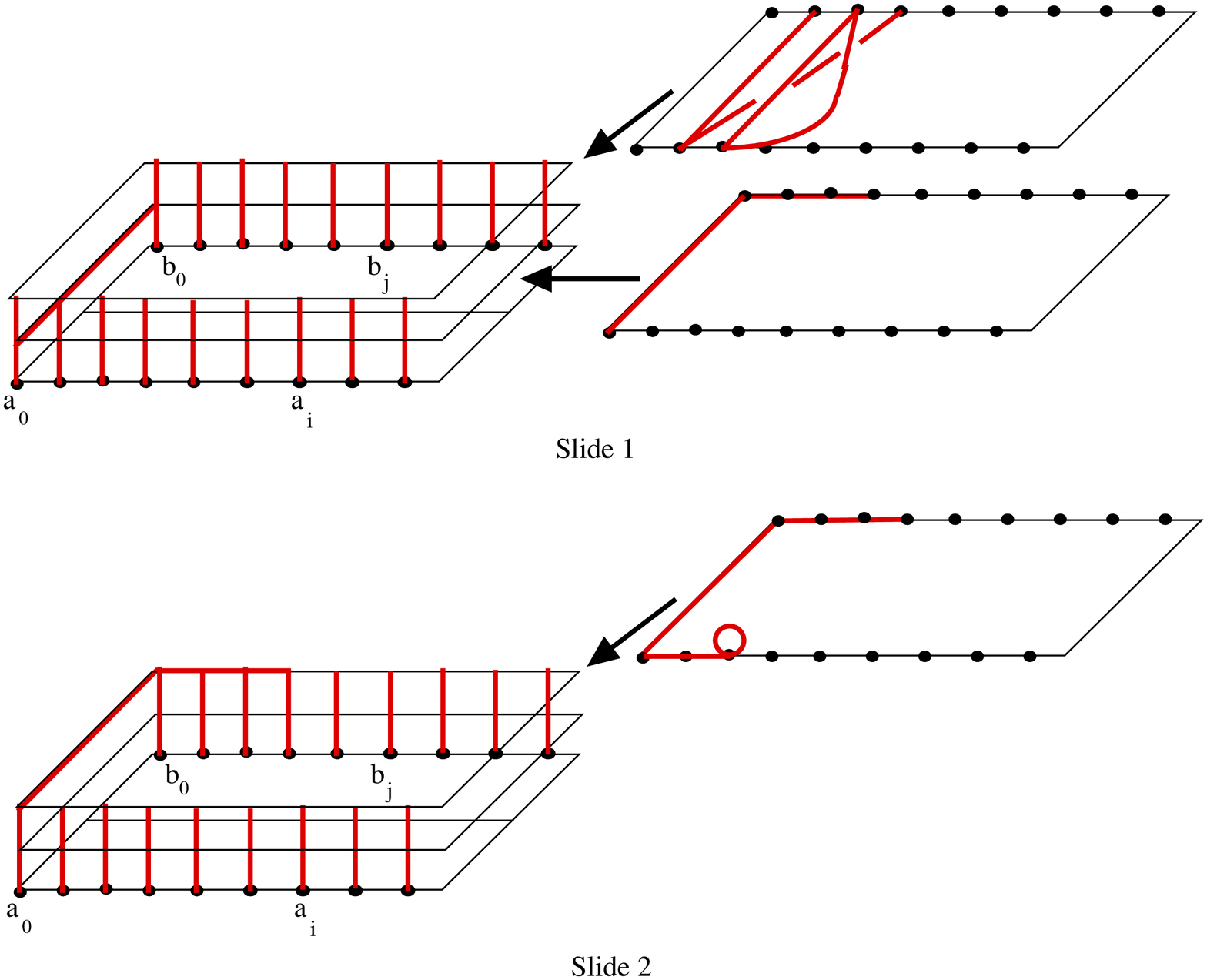}
\caption{Sliding $\Lll$ in the cube by layer} \label{fig:bipartite2}
\end{figure}

\begin{proof} We will assume $\Lll$ is connected; if not, the 
following argument can be carried out in each component separately.

Place the designated edge $e$ as described.  Denote its ends 
by $a_{0} = (0, -1, 0)$ and $b_{0} = (0, 1, 0)$.  The $\Lll$-distance 
between two vertices in $\Lll$ will mean the number of edges in the 
shortest path between them.  With no loss, order the indices 
of the remaining vertices $a_{i}, 1 \leq i \leq |A|-1$ of $A$ 
subordinate to their $\Lll$-distance from $a_{0}$, i.  e.  so that, 
for any pair of indices $i_{1}$ and $i_{2}$, if $a_{i_{1}}$ is closer 
in $\Lll$ to $a_{0}$ than $a_{i_{2}}$ is, then $i_{1} < i_{2}$.  (We 
do not care how vertices are ordered among those that are $\Lll$-equidistant 
from $a_{0}$.)  Similarly order the indices of the remaining vertices 
$b_{j}, 1 \leq j \leq |B|-1$ of $B$ subordinate to their 
$\Lll$-distance from $a_{0}$.  After this reordering, place each 
$a_{i}$ at the point $(i/|A|, -1, 0)$ and each $b_{j}$ at the point 
$(j/|B|, 1, 0)$.

At each vertex of $\Lll$ add a vertical (i.  e.  $z$-parallel) arc of 
length $1$.  That is, attach to each $(i/|A|, -1, 0)$ the arc $\{ 
(i/|A|, -1) \} \times [0, 1]$ and to each $(j/|B|, 1, 0)$ the arc $\{ 
(j/|B|, 1) \} \times [0, 1]$.  In order to simplify somewhat the 
description of $\Lll$, the edges of $\Lll$ will originally be placed 
so that they are horizontal (i.  e.  parallel to the $x-y$ plane) with 
ends on these vertical arcs.  $\Lll$ is then finally recovered from 
the simplified description by collapsing the vertical arcs $\{ (i/|A|, 
-1) \} \times [0, 1]$ and $\{ (j/|B|, 1) \} \times [0, 1]$ back down to 
$A$ and $B$ respectively.

Let $\ell$ be the maximal $\Lll$-distance of any vertex in $\Lll$ from
$a_{0}$.  We will place the edges of $\Lll$ in a sequence of $\ell$
stages; the edges placed at the $k^{th}$ stage lie near the horizontal
square $I \times [-1, 1] \times \{ k/\ell \}$.  Specifically, at the
$k^{th}$ stage select all edges of $\Lll$ which have the property that
their most $\Lll$-distant end is a $\Lll$-distance $k$ from $a_{0}$. 
(The other end of each selected edge must then be $\Lll$-distance
$k-1$ from $a_{0}$.)  If there are $p$ such edges, select a sequence
of $p$ horizontal planes whose height (i.  e. $z$-coordinate) is near
$k/\ell$ and place each edge in a separate horizontal plane, as a
linear edge connecting the appropriate $a_{i}$ to the appropriate
$b_{j}$, with parallel edges on adjacent horizontal planes.  The
linear embedding ensures that each edge is monotonic in the
$y$-coordinate, a fact that is unchanged when the vertical arcs $\{
(i/|A|, -1) \times [0, 1]$ and $(j/|B|, 1) \} \times [0, 1]$ are
collapsed to $A$ and $B$ to create $\Lll$.  We have thereby described
an embedding of $\Lll$ into the cube that clearly satisfies the first
three requirements.  See Figure \ref{fig:bipartite}

\begin{figure}[tbh] 
\centering
\includegraphics[width=1.0\textwidth]{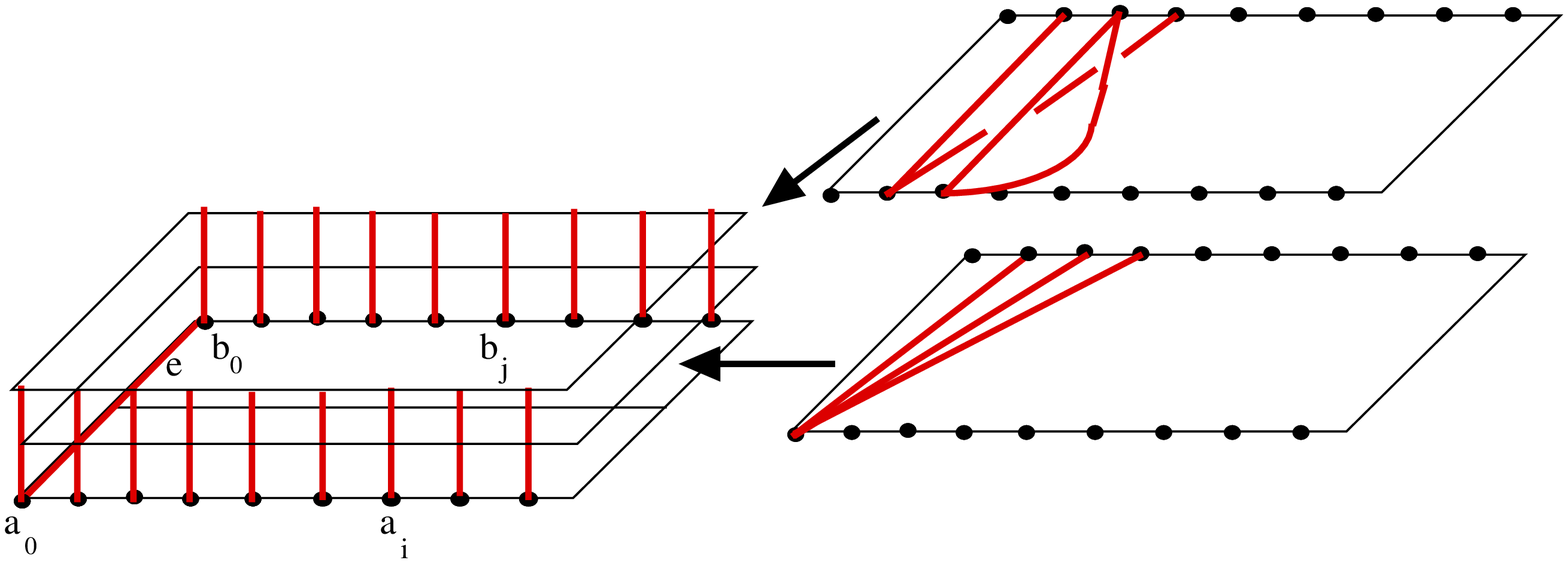}
\caption{Putting $\Lll$ in the cube in layers} \label{fig:bipartite}
\end{figure}

It remains to describe how the edges of $\Lll$ can be slid and 
isotoped, without moving the vertices $A, B$ or the edge $e$, so that 
afterwards the resulting graph lies entirely in the face $I \times 
[-1, 1] \times \{ 0 \}$.  The description of this sliding mimics the 
$k$ stages of the construction of $\Lll$ and we will describe them in 
the graph above as if we had not collapsed the vertical arcs, but also 
mostly focusing on the $x-y$ coordinates.  

At the first stage of the construction above, exactly those edges with
one end on $a_{0}$ are added, near the horizontal plane $z = 1/\ell$. 
By our choice of ordering of the $b_{j}$, the other ends of these
edges lie exactly on the vertices $b_{0}, \ldots, b_{q}$, for some $q
\geq 0$.  (If any two of these edges are parallel, slide one over the
other to form a tiny circle which we may henceforth ignore).  Then, if
$q > 0$ the rightmost edge, i.  e. that connecting $a_0$ to $b_{q} =
(q/|B|, 1)$ may be slid over the edge connecting $a_0$ to $b_{q-1}$
until instead it is just the straight interval between $b_{q-1}$ and
$b_{q}$, i.  e. the interval $[q - 1, q] \times \{ 1 \}$.  Continue in
this manner until all the edges but $e$ have been slid to the line $y
= 1$ to constitute the single interval $[0, q] \times \{ 1 \}$, still
in the plane $z = 1/\ell$.  See Figure \ref{fig:bipartite2}.  Now
slide all these edges up vertically to height just below $z = 2/\ell$
and begin the second stage.

Because of our ordering of the $a_{i}$, there is a $p \geq 1$ so that
the vertices $a_{1}, \ldots, a_{p}$ constitute exactly the ends in $A$
of edges included at the second stage.  Moreover the other end of each
such edge lies among the $b_{0}, \ldots, b_{q}$ which, after the
slides we have done on the edges of the first stage, all lie on the
$L$-shaped graph $e \cup ([0, q] \times \{ 1 \})$.  This $L$-shaped
graph gives a way, much as above, of sliding the edges added at the
second stage until they are either tiny circles (henceforth ignored)
or constitute the straight line from $a_{0}$ to $a_{p}$, i.  e. the
line $[0, p] \times \{ -1 \} \times 2/\ell$.  See Figure
\ref{fig:bipartite2}.  Now slide this whole graph vertically up until
it is near the plane $z = 3/\ell$ and continue the process.  By the
time we have reached the $\ell^{th}$ stage, the graph consists (now at
height $z = 1$) of arcs in the lines $y = \pm 1$ that contain all the
vertices, together with the original edge $e$ between $a_{0}$ and
$b_{0}$ (and some tiny circles), all of which then lie in the square
$I \times [-1, 1] \times \{ 1 \}$.  Now collapse the vertical
direction, bringing the graph down to $I \times [-1, 1] \times \{ 0
\}$.  This process (when reinterpreted as slides on the actual
embedding of $\Lll$, in which the vertical arcs do not appear)
verifies the last numbered condition.
\end{proof}

\section{Braid equivalence and unknotted graphs} \label{section:braid}

Suppose $(M, h)$ is a planar presentation and $t \in R$ is a regular 
value of $h$.  Cut $M$ along $P^{t}$ and reattach the two copies of 
$P^{t}$ by an orientation preserving homeomorphism $P^{t} \map P^{t}$ 
that is the identity on the circles $\bdd P^{t}$.  The result is a 
possibly new manifold $M'$ and a planar presentation $h': M' \map R$.  
Note that $h'|\bdd M' = h|\bdd M$.  The two planar presentations $(M, 
h)$ and $(M', h')$ are called {\em braid equivalent}.  More generally, 
two planar presentations $(M, h)$ and $(M', h')$ are called braid 
equivalent if one is obtained from the other by a finite sequence of 
such operations, called braid moves.

Under such braid moves, many more $3$-manifolds with 
planar presentation can be made unknotted graph complements.  The 
following lemma illustrates why.  The setting is this: Suppose $N_{A}$ 
and $N_{B}$ are each homeomorphic to either $B^{3}$ or $S^{2} \times 
I$ and $P_{A}$ (resp.  $P_{B}$) is a sphere component of the boundary 
of $N_{A}$ (resp $N_{B})$.  Let $N$ be obtained from identifying 
$P_{A}$ with $P_{B}$ (so in particular $N$ is also homeomorphic to 
either $B^{3}$ or $S^{2} \times I$). Suppose 
further that $\Lll \subset N$ is a properly embedded graph that is in 
general position with respect to $P$; let $\Lll_{A} = \Lll 
\cap N_{A}$ and $\Lll_{B} = \Lll \cap N_{B}$.

\begin{lemma} \label{lemma:sum} If both $\Lll_{A}$ and $\Lll_{B}$ are 
unknotted graphs, then there is a homeomorphism $\phi: P_{A} \map 
P_{B}$ such that 

\begin{enumerate}
    
    \item $\phi$ is the identity on the points $\Lll \cap 
P_{A}$ and 

\item $\Lll_{A} \cup \Lll_{B}$ is an unknotted graph in $N_{A} 
\cup_{\phi} N_{B}$.
\end{enumerate}
\end{lemma}

\begin{proof} The case in which both $N_{A}$ and $N_{B}$ are copies 
of $S^{2} \times I$ is representative (and in fact the most 
difficult) and it will be convenient to take $N_{A} = S^{2} \times 
[-2, 0]$ and $N_{B} = S^{2} \times [0, 2]$.

Construct an abstract bipartite graph $G$ with vertex sets $A$ and $B$ 
as follows: There is a vertex in $A$ (resp $B$) for every component of 
$\Lll_{A}$ (resp $\Lll_{B}$).  There is an edge for every point $c$ in 
$P_{A} \cap \Lll = P_{B} \cap \Lll$.  Identify the ends of the edge 
corresponding to $c$ to the points in $A$ and $B$ corresponding to the 
components in $\Lll_{A}$ and $\Lll_{B}$ on which $c$ lies.  Imbed $G$ 
in the cube $I \times [-1, 1] \times I$ as described in Lemma 
\ref{lemma:bipartite} and embed the cube in $S^{2} \times [-2, 2]$ 
with the $x-z$ square cross-section of the cube lying in the $S^{2}$ 
factor and the $y$-coordinate of the cube projecting to the interval 
factor in $S^{2} \times [-1, 1] \subset S^{2} \times [-2, 2]$.

For each vertex $v$ in $\bdd \Lll \cap (S^{2} \times \{ -2 \})$ add a 
monotone edge $e_{v} \subset S^{2} \times [-2, -1]$ to $G$ with one 
end of $e_{v}$ on $v$ and the other end on the vertex in $A$ 
corresponding to the component of $\Lll_{A}$ on which $v$ lies.  
Similarly, add a monotone edge in $S^{2} \times [1, 2]$ for each 
vertex in $\bdd \Lll \cap (S^{2} \times \{ 2 \})$, with one edge on 
the vertex and the other on the appropriate vertex in $B$.  Call the 
resulting graph $G_{+}$.  See Figure \ref{fig:braidmove}.

\begin{figure}[tbh] 
\centering
\includegraphics[width=0.5\textwidth]{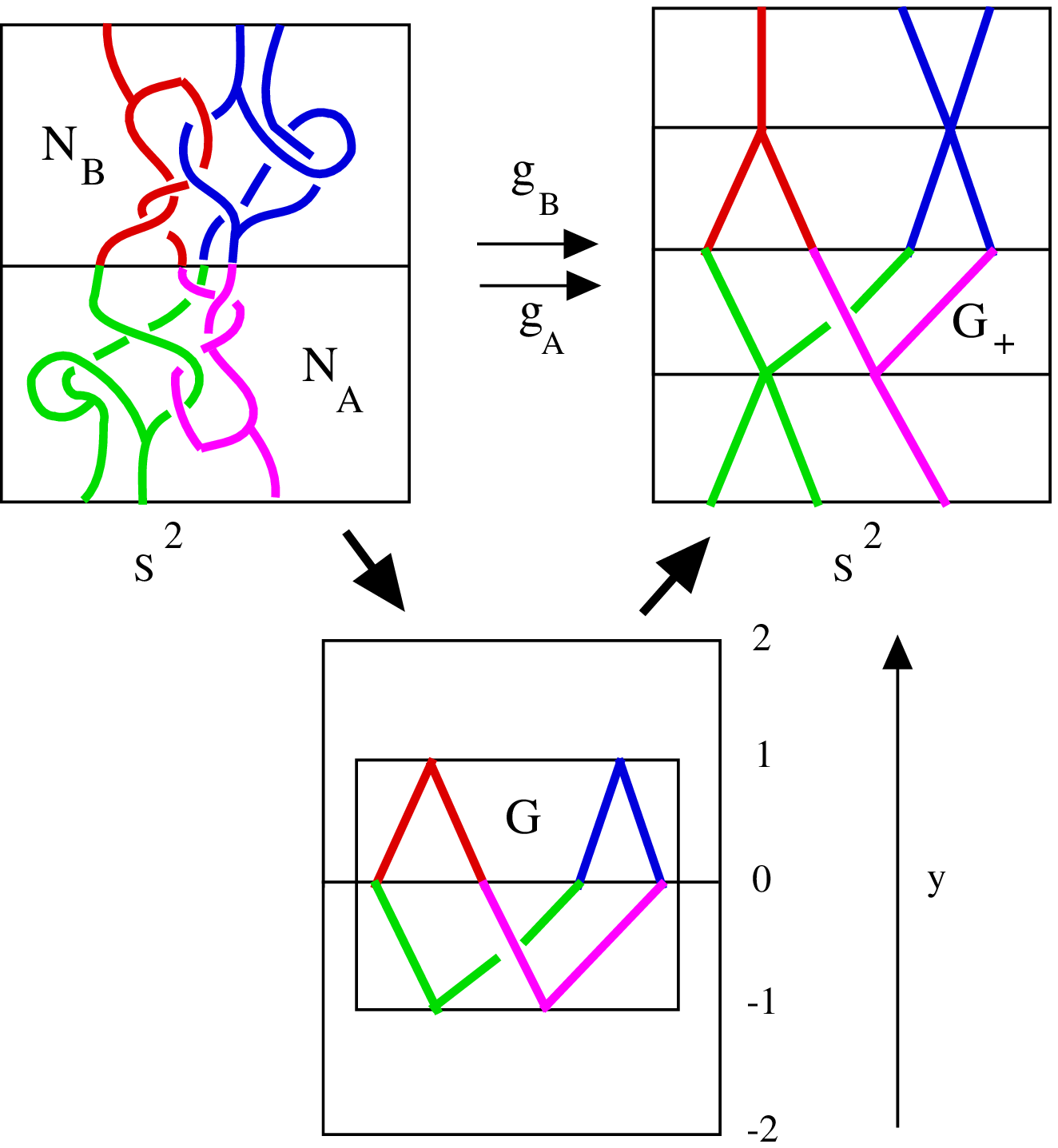}
\caption{A braid move makes $\Lll$ unknotted} \label{fig:braidmove}
\end{figure} 

The graph $G_{+}$ as embedded, has three important properties:

\begin{itemize}
    
    \item It follows from Lemma \ref{lemma:bipartite} that $G_{+}$ is an 
    unknotted graph in $N$.
    
    \item It follows from Proposition \ref{prop:unique} and Example 
    \ref{exam:concave} that, perhaps after adding some tiny 
    circles to $G_{+}$, $G_{+} \cap N_{A}$ and $\Lll_{A}$ are equivalent 
    unknotted graphs, via an equivalence that is the 
    identity near $\bdd \Lll_{A}$.
    
    \item Similarly, perhaps after adding some tiny 
    circles to $G_{+}$, $G_{+} \cap N_{B}$ and $\Lll_{B}$ are equivalent 
    unknotted graphs via an equivalence that is the 
    identity near $\bdd \Lll_{B}$.
    
    \end{itemize}
    
    Let $g_{A}: P_{A} \map P_{A}$ and $g_{B}: P_{B} \map P_{B}$ be the 
    homeomorphisms given by the latter two equivalences.  Let $\phi = 
    g_{B} g_{A}^{-1}: P_{A} \map P_{B}$.  Then the construction $N_{A} 
    \cup_{\phi} N_{B}$ changes $\Lll$ to a graph equivalent to $G_{+}$, 
    which is unknotted.
\end{proof}

This has as an immediate corollary, analogous to Example 
\ref{exam:pileon2}.  Suppose $(M, h)$ is a planar presentation of a 
manifold $M$ and $M^{[a,c]}$ is a component of $h^{-1}([a, c])$.  We 
again will use the following notation: for $J$ a subinterval of $[a, 
c]$ let $M^{J} = M^{[a,c]} \cap h^{-1}(J)$ and for $t \in [a, c]$ let 
$Q^{t} = M^{[a,c]} \cap P^{t}$.

\begin{cor} \label{cor:medexam} Suppose $(M, h)$ is a planar 
presentation of a manifold $M$ and $M^{[a,c]}$ is a 
component of $h^{-1}([a, c])$.  

Suppose that for some $b \in [a, c]$

\begin{enumerate}

     \item $Q^{b}$ is connected
    
    \item the pair $(M^{[b,c]}, Q^{b} \cup Q^{c})$ is an unknotted 
    graph complement.
    
    \item  all saddles in  $M^{[a, b]}$ are nested. 

\end{enumerate}

Then $M^{[a, c]}$ is braid-equivalent to an unknotted graph complement.
\end{cor}

\begin{proof} The proof is by induction on the number of critical 
points of $h$ on $\bdd M$ that occur in $M^{[a, b]}$.  If there are 
none then of course $M^{[a, c]} \cong M^{[b, c]}$ and there is nothing 
to prove.  If the highest singularity in $M^{[a, b]}$ is a maximum or 
a minimum (necessarily an internal max or min since $Q^{b}$ is 
connected and all saddles in $M^{[a, b]}$ are nested) then for $t$ 
just below the corresponding critical value, $M^{[t, c]}$ is a 
standard graph complement and we are done by induction.  Similarly, if 
the highest critical value in $[a, b]$ is a (nested) upper saddle then 
apply Example \ref{exam:pileon2} to complete the inductive step.  

The only remaining case is when the highest critical point is a lower 
saddle, i.  e.  it suffices to consider the case in which the only 
critical point in $M^{[a, b]}$ is a single nested lower saddle.  But 
even in the more general case that all saddles in $M^{[a, b]}$ are 
nested lower saddles, the proof is an immediate consequence of Lemma 
\ref{lemma:sum} and Example \ref{exam:medexam} with the latter applied 
to $M^{[a, b]}$, which has no upper saddles.
\end{proof}

\bigskip

\begin{figure}[tbh] 
\centering
\includegraphics[width=.3\textwidth]{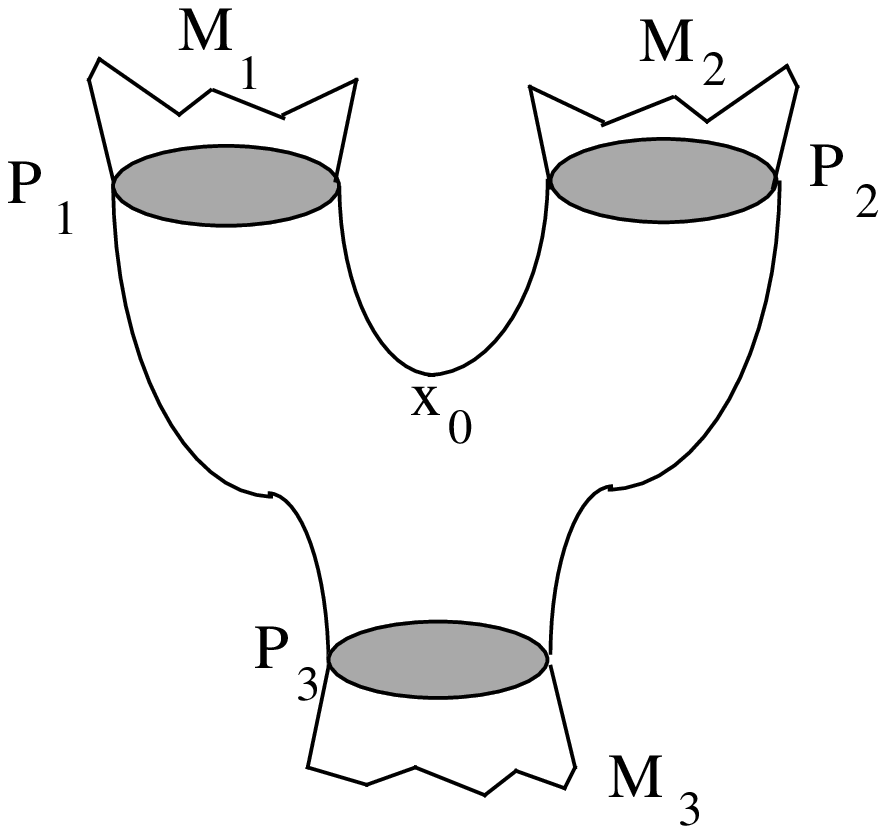}
\caption{} \label{fig:saddle2}
\end{figure}

We hope next to understand what happens to planar presentations of
unknotted graph complements at unnested saddles.  So let $a$ be a
critical value with corresponding critical point $x_{0}$, an unnested
saddle.  For small $\eee$, let $M^{[a - \eee, a + \eee]}$ be the
component of $h^{-1}([a - \eee, a + \eee])$ that contains $x_{0}$. 
Then for, say, a lower saddle, $P^{a + \eee}$ intersects $M^{[a -
\eee, a + \eee]}$ in two connected planar surfaces denoted $P_{1}$ and
$P_{2}$ and $P^{a - \eee}$ intersects $M^{[a - \eee, a + \eee]}$ in a
single connected planar surface $P_{3}$.  The roles of $\pm \eee$ are
reversed for an upper saddle.  We will be interested only in the case
in which each $P_{i}$ separates $M$.  The component of $M - P_{i}$ not
containing the saddle point will be denoted $M_{i}$.  See Figure
\ref{fig:saddle2}.

One situation is easy to understand:

\begin{lemma} \label{lemma:straight} If $(M_{1}, P_{1})$ and $(M_{2}, 
P_{2})$ are unknotted graph complements in $B^{3}$, then so is $(M - 
interior(M_{3}), P_{3})$, the component of the complement of $P_{3}$ 
that contains both $M_{1}$ and $M_{2}$.  
\end{lemma}

\begin{proof} For this proof, a useful model of an unknotted graph in 
$B^{3}$ is this: In a cube $I \times I \times I$, let $\Lll$ be a 
subgraph of the square $I \times I \times \{ 1/2 \}$ with a single 
boundary vertex on the top $I \times \{ 1 \} \times \{ 1/2 \}$ and the 
rest on the bottom $I \times \{ 0 \} \times \{ 1/2 \}$. (Here
projection to the $y$-coordinate models height $h$). Then the 
complement of $\Lll$ in the cube is an unknotted graph complement with 
planar part on the bottom of the box, namely $(I \times \{ 0 \} \times 
I ) - \eta(\bdd \Lll)$.  Note that the sides of the box are not in the 
planar part but rather it is much as if the vertex of $\bdd \Lll$ at the 
top of the box is stretched over the top and all sides of the box.  
See Figure \ref{fig:straight}.

\begin{figure}[tbh] 
\centering
\includegraphics[width=.8\textwidth]{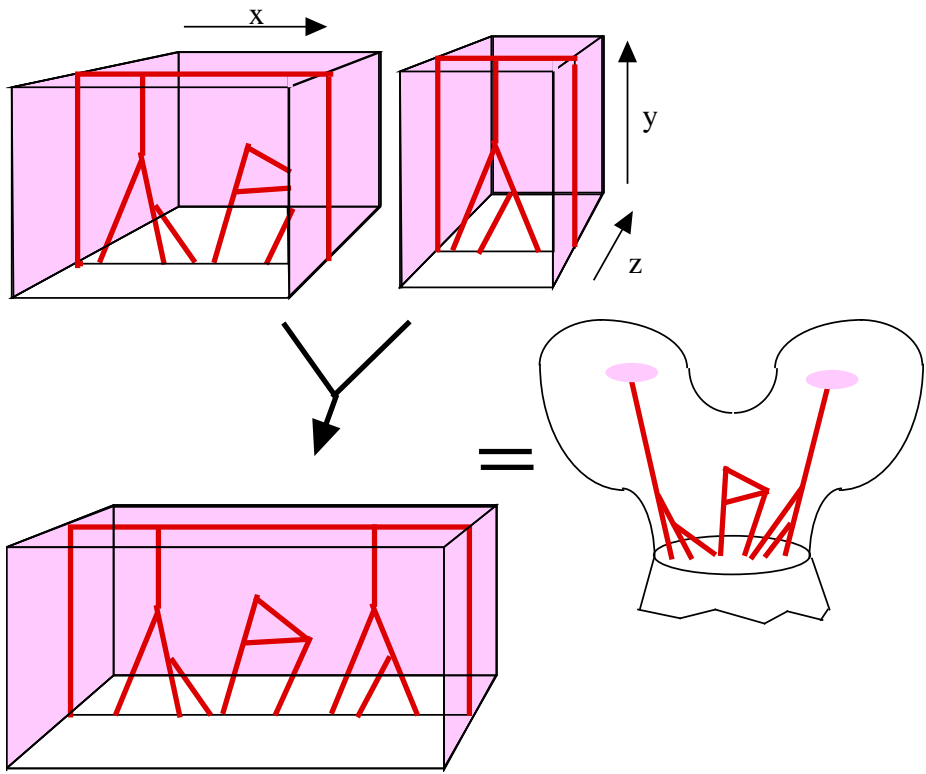}
\caption{Modelling an unnested lower saddle} \label{fig:straight}
\end{figure}

The effect of passing through an unnested lower saddle is to take two such 
boxes (each containing one $P_{i}$ on its bottom) and glue the side 
$\{ 1 \} \times I \times I$ of one to the side $\{ 0 \} \times I 
\times I$ of the other, obtaining a graph complement with planar part 
the boundary sum of the original two planar parts.  The result is 
again a cube with the same sort of graph deleted, with the sole 
difference that now there are two boundary vertices of the graph on 
the top of the box.  But since the top of the box is entirely disjoint 
from the planar part of the graph complement, up to graph complement 
equivalence, nothing is changed by sliding one top boundary vertex to 
the other along the top arc $I \times \{ (1, 1/2) \}$, and then 
sliding an end of one edge down the end of the other, after which 
there is again a single boundary vertex on the top.  In particular, 
the result is again an unknotted graph complement in the cube.
\end{proof}
    
A much harder situation to analyze is this:

\begin{lemma} \label{lemma:turn} If $(M_{3}, P_{3})$ and $(M_{2},
P_{2})$ are unknotted graph complements, then the pair $(M -
interior(M_{1}), P_{1})$ is braid equivalent to an unknotted graph
complement.
\end{lemma}

\begin{proof}  The initial difficulty is to determine a good model 
for what we are trying to show, analogous to the model in Lemma 
\ref{lemma:straight}.  Let $\Lll_{A}, \Lll_{B}$ be unknotted graphs in 
the $3$-ball whose complements give $M_{2}$ and $M_{3}$ respectively.  
Inspired by the model above (with the $y$-coordinate again modelling 
the height function $h$, but this time for an unnested upper saddle) 
choose two cubes $C_{2}, C_{3}$ in $R^{3}$, as follows (see Figure 
\ref{fig:turn}:

\begin{enumerate}
    
    \item $C_{2} = [0, 1] \times [-2, 0] \times [-1, 1]$
    
    \item $C_{3} = [0, 2] \times [0,2] \times [-1, 1]$.
      
\end{enumerate}
    
    Let $C_{\cup} = C_{2} \cup C_{3}$, which is itself homeomorphic to 
    a $3$-ball.
    
Construct an abstract bipartite graph $G$ with vertex sets $A$ and $B$ 
as follows: There is a vertex in $A$ (resp $B$) for every component of 
$\Lll_{A}$ (resp $\Lll_{B}$).  There is an edge for every component 
$c$ of $\bdd P_{2}$.  Identify the ends of such an edge to the points 
in $A$ and $B$ that represent the components of $\Lll_{A}$ and 
$\Lll_{B}$ on which $c$ lies.  Imbed $G$ in the cube $I \times [-1, 1] 
\times I \subset C_{\cup}$ as described in Lemma \ref{lemma:bipartite} 
with the special edge $e$ chosen to be that which corresponds to the 
boundary component of $P_{2}$ which is incident to the saddle 
singularity.  (Notice that $e$ lies on the face $\{ 0 \} \times I 
\times I$ of $\bdd C_{\cup}$.)

The vertices of $B$ are strung out along the interval $I \times \{ 1
\}$ in the $x-y$ plane, with all but the vertex of $e$ lying in the
interor of $C_{\cup}$.  Add edges to $G$ that connect these vertices
of $B$ linearly to the corresponding vertices in the interval $\{ 1 \}
\times [1, 2]$ in the $x-y$ plane.  Explicitly, add an edge $e_{j}, j
= 1,\ldots, |B| - 1$ that connects the point $b_{j} = (j/|B|, 1, 0)
\in B$ to the point $(1, 2 - j/|B|, 0)$.  Next add edges that connect
these points linearly to a collection $B'$ of points in the line $[1,
2] \times \{ (0, 0) \} \subset \bdd C_{3}$.  This collection $B'$ is
chosen so that each point in $B'$ corresponds to a boundary component
of $P_{1}$, other than the one containing the saddle singularity. 
Equivalently, each point $b' \in B'$ corresponds to a vertex in $\bdd
\Lll_{B}$ that doesn't also naturally correspond to a vertex in $\bdd
\Lll_{A}$.  Such a boundary vertex lies on a component of $\Lll_{B}$
to which a vertex $b_{j}$ has been assigned; append a linear edge in
$[1, 2] \times [0,2] \times [-1, 1]$from the {\em other} end of
$e_{j}$ to $b'$.  (We pick the ordering of $B'$ in the interval $[1,
2] \times \{ (0, 0) \}$ so that these edges do not intersect.) 
Finally, append an appropriate number of tiny circles to $G \cap
C_{2}$ and $G \cap C_{3}$ so that each component has the same Euler
number as the corresponding component of $\Lll_{A}$ and $\Lll_{B}$. 
Let $G_{+}$ be the graph in $C_{\cup}$ given by this construction. 
Note that it is a proper graph in $C_{\cup}$ whose planar part
$P_{\cup}$ we take to be $([1,2] \times \{0 \} \times [-1, 1]) -
\eta(G_{+}$), i.  e. the complement of $G_{+} \cup C_{2}$ in the
bottom face of $C_{3}$.  See Figure \ref{fig:turn}.

\begin{figure}[tbh]
\centering
\includegraphics[width=.7\textwidth]{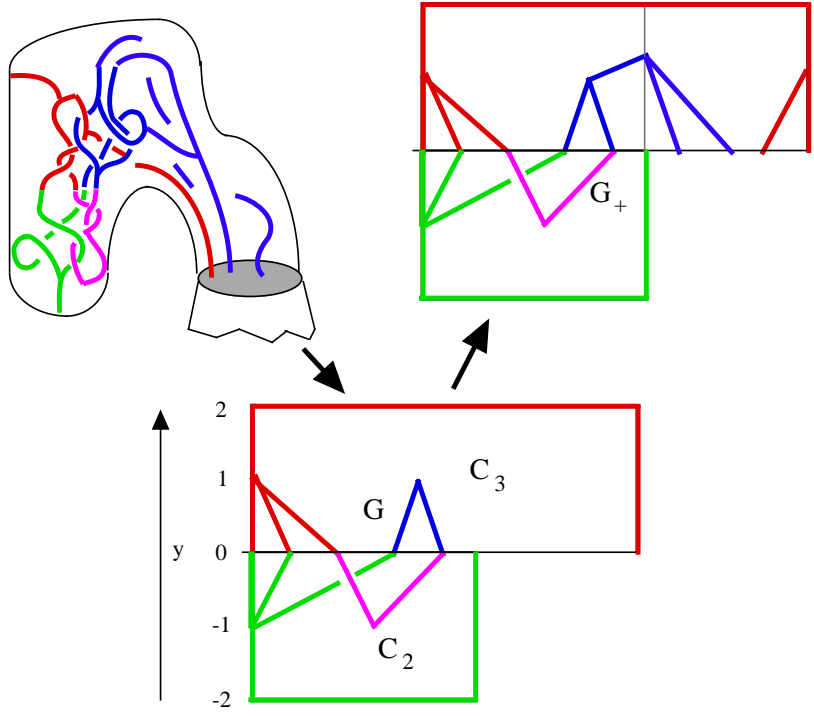}
\caption{Modelling a complicated unnested upper saddle} \label{fig:turn}
\end{figure} 

Let $P_{x-z} \subset R^{3}$ denote the plane $y = 0$.  The graph
$G_{+}$ has been constructed to have these properties:

\begin{enumerate}
    
    \item For $i = 2, 3$, the graph $\Ggg_{i} = G_{+} \cap C_{i}$ is 
    unknotted, with planar part $(P_{x-z} \cap C_{i}) - \eta(G_{+})$, 
    by Example \ref{exam:smallexam}.
    
     \item Each component of $\Ggg_{2}$ (resp $\Ggg_{3}$) is 
     homeomorphic to a corresponding component of $\Lll_{A}$ (resp 
     $\Lll_{B}$) so that the homeomorphisms agree, where they are 
     simultaneously defined, namely on $\bdd \Ggg_{2} \subset 
     \bdd \Ggg_{3}$.
    
     \item The graph $G_{+} \subset C_{\cup}$ is unknotted by Lemma 
     \ref{lemma:bipartite}.
    
\end{enumerate}
    
The first two properties guarantee (via Proposition \ref{prop:unique}) 
that there is a homeomorphism of pairs $(M_{i}, P_{i}) \cong (C_{i} - 
\eta(G_{+}), (P_{x-z} \cap C_{i}) - \eta(G_{+}))$.  In particular, much 
as in Lemma \ref{lemma:sum}, $M_{3}$ can be cut off from $M$ 
and reattached so that the pair $(M - interior(M_{1}), P_{1})$ becomes 
pairwise homeomorphic to $(C_{\cup} - \eta(G_{+}), P_{\cup})$.  But 
since $G_{+}$ is unknotted, the latter is an unknotted graph 
complement.  Hence $(M - interior(M_{1}), P_{1})$ is braid equivalent to a 
standard graph complement.
\end{proof}

\section{Heegaard reimbedding} \label{section:heegaard}

\begin{theorem} \label{thm:main} Suppose $(M, h)$ is a planar 
presentation of a $3$-manifold with connectivity graph a tree.  Then 
$(M, h)$ is braid-equivalent to an unknotted graph complement.
\end{theorem}

\begin{proof} If the connectivity graph $\Ggg$ is a vertex (i.  e.  
all saddles are nested) the result follows easily from Corollary 
\ref{cor:medexam}.  So we will assume that $\Ggg$ has at least one 
edge.  In that case, Lemma \ref{lemma:sum} demonstrates that the proof 
of the theorem will follow from the proof of the following relative 
version.
\end{proof}

\begin{prop} Suppose $(M, h)$ is a planar presentation of a
$3$-manifold and $\Ggg$ is its connectivity graph.  Suppose $\ggg
\subset \Ggg$ is an edge such that a component $\Ggg_{0} \subset \Ggg$
of the complement of $\ggg$ is a tree.  Let $P_{\ggg} \subset M$ be
the planar surface corresponding to $\ggg$ and $M_{0} \subset M$ be
the component of $M - P_{\ggg}$ that corresponds to $\Ggg_{0}$.  Then
$(M_{0}, P_{\ggg})$ is braid equivalent to an unknotted graph
complement.
\end{prop}

\begin{proof} The proof will be by induction on the number of edges in 
$\Ggg_{0}$.  Let $v$ be the vertex of $\Ggg_{0}$ that is incident to 
$\ggg$ and, in the terminology of Lemma \ref{lemma:vertex}, let 
$M_{v}$ be the component of $M - \cup_{i=1}^{n} P^{s_{i}}$ 
corresponding to $v$, with $h(M_{v}) = [s_{i}, s_{i+1}]$.  We will 
assume that the unnested saddles at heights $s_{i}$ and $s_{i+1}$ both 
involve the particular component $M_{v}$, since the argument is 
easier if either or both do not.  

Without loss of generality we will assume that the planar surface
corresponding to the edge $\ggg$ is at the bottom of $M_{v}$, i.  e.
near height $s_{i}$.  Consider first the saddle $x_{+}$ at height
$s_{i+1}$.  Let $Q$ be the connected planar surface $M_{v} \cap
P^{(s_{i+1} - \eee)}$ and $M_{Q}$ be the component of $M - Q$ that
contains $x_{+}$.  If $x_{+}$ is an upper (unnested) saddle then $Q$
corresponds to an edge in $\Ggg_{0}$ and so by inductive assumption
the pair $(M_{Q}, Q)$ is braid equivalent to an unknotted graph
complement.  See Figure \ref{fig:tree}a.  If $x_{+}$ is a lower saddle then
the two contiguous components of $h^{-1}(s_{i+1} + \eee)$ each
represent edges in $\Ggg_{0}$ and $(M_{Q}, Q)$ is again an unknotted
graph complement by inductive assumption combined with Lemma
\ref{lemma:straight}.  See Figure \ref{fig:tree}b.  So in any case, $(M_{Q},
Q)$ is braid equivalent to an unknotted graph complement.

\begin{figure}[tbh] 
\centering
\includegraphics[width=.8\textwidth]{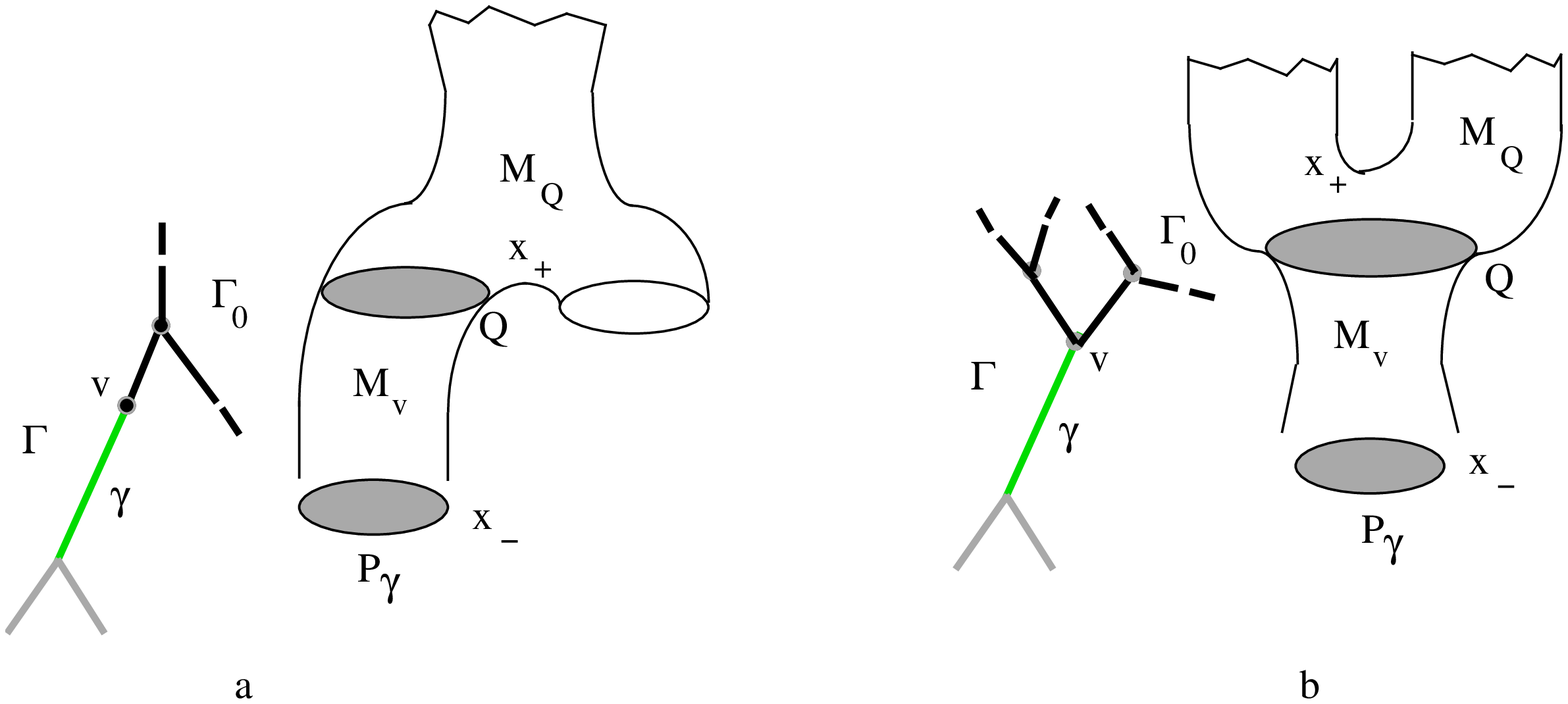}
\caption{} \label{fig:tree}
\end{figure}

Now consider the saddle $x_{-}$ at height $s_{i}$. See Figure 
\ref{fig:tree2}. If it's a lower
saddle, then the planar surface $P^{s_{i} + \eee} \cap M_{v}$ is
$P_{\ggg}$, the planar surface corresponding to the edge $\ggg$ and
the proposition follows from Corollary \ref{cor:medexam}.  If the
saddle $x_{-}$ is an upper saddle, then $P_{\ggg}$ is one of the two
connected planar surfaces in $P^{s_{i} - \eee}$ contiguous to the
saddle.  Let $P_{\ggg'}$ be the other one, with corresponding edge
$\ggg' \subset \Ggg_{0}$, and let $P$ be the connected planar surface
$M_{v} \cap P^{s_{i} + \eee}$.  Now by inductive assumption, the
component of $M_{0} - P_{\ggg'}$ not containing $x_{-}$ is an unknotted
graph complement and by Corollary \ref{cor:medexam} so is the
component of $M_{0} - P$ not containing $x_{-}$.  Then the proposition
follows from Lemma \ref{lemma:turn}.

\begin{figure}[tbh] 
\centering
\includegraphics[width=.8\textwidth]{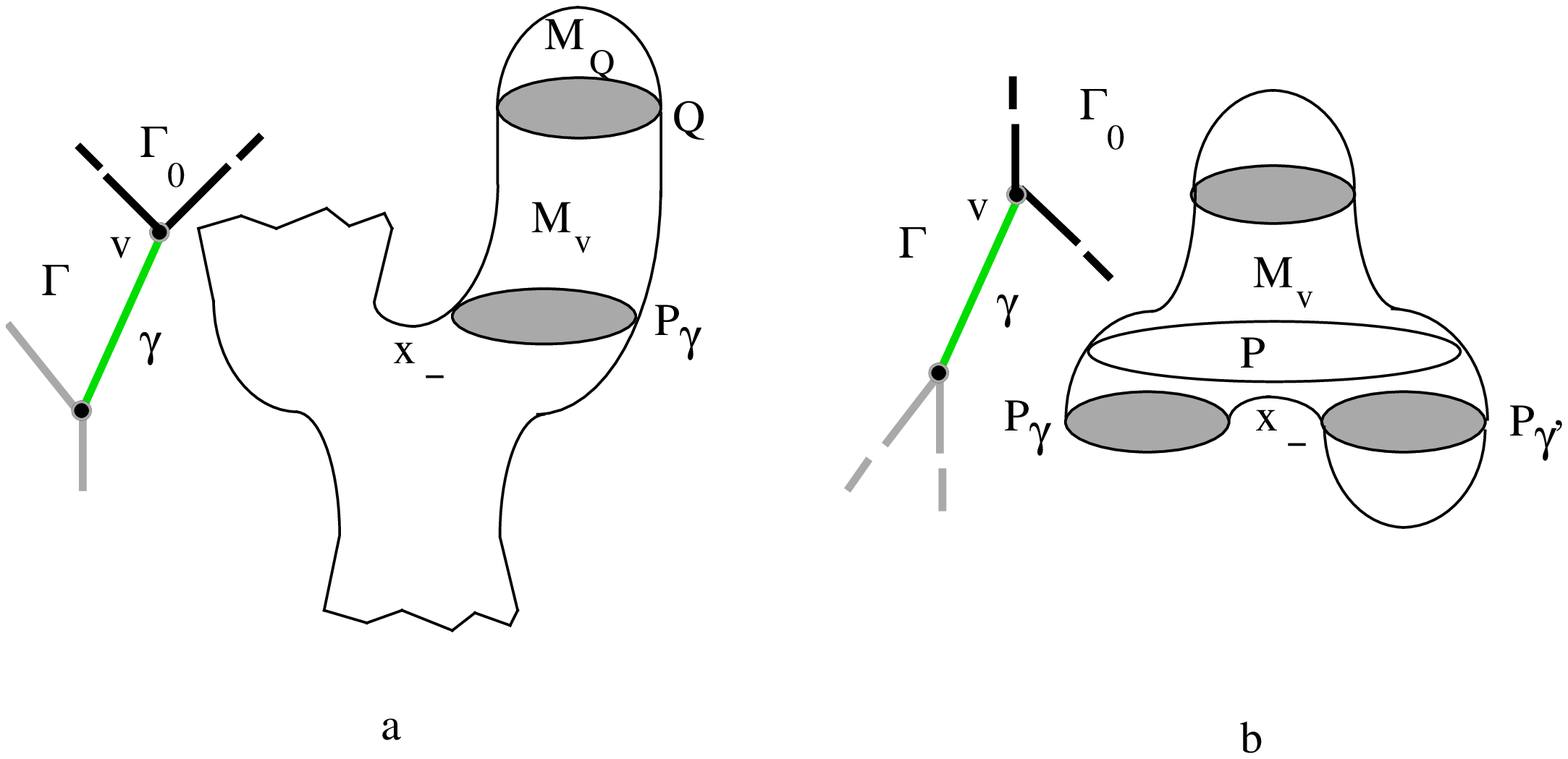}
\caption{} \label{fig:tree2}
\end{figure}
\end{proof}

\begin{cor} \label{cor:reimbed} Suppose $N \subset S^{3}$, $p: S^{3} 
\map R$ is the standard height function, $N$ contains both poles, and 
the connectivity graph of $S^{3} - N$ is a tree (so in particular 
$S^{3} - N$ is connected).  Then there is an embedding $f: N \map 
S^{3}$ so that

\begin{enumerate} 

\item $p = p f$ on $N$, i. e. $f$ preserves height and

\item $S^{3} - f(N)$ is a connected sum of handlebodies.

\end{enumerate}
\end{cor}

\begin{proof} It follows from Theorem \ref{thm:main} that $M = S^{3} - 
N$ is braid-equivalent to a connected sum of handlebodies.  We will 
show that a braid move on $M$ defines a reimbedding of 
$N$.

Let $S^{t}$ be the $2$-sphere $p^{-1}(t)$ and $P^{t} = S^{t} - N =
S^{t} \cap M$.  Then a braid move of $M$ at a generic level $t$ is
given by cutting $M$ open along $P^{t}$ and then reattaching $P^{t}$
to itself by a homeomorphism $\phi: P^{t} \map P^{t}$ that is the
identity on $\bdd P^{t}$.  In particular, the homeomorphism $\phi$
extends via the identity on $S^{2} - P^{t}$ to a self-homeomorphism of
$S^{2}$.  But any (orientation preserving) self-homeomorphism of the
sphere is isotopic to the identity, so in fact there is a
level-preserving self-homeomorphism $S^{2} \times [t - \eee, t+ \eee]$
that is the identity on one end and the extended $\phi$ on the other. 
Use this self-homeomorphism to redefine the embedding of $N$ in the
region $h^{-1}[t - \eee, t+ \eee]$.  The effect on the complement $M$
is to do the original braid move.
\end{proof}

\begin{cor} \label{cor:heegaard} Suppose $p: S^{3} \map R$ is the 
standard height function and $H \subset S^{3}$ is a handlebody for 
which horizontal circles constitute a complete collection of meridian 
disk boundaries.  Then there is a reimbedding $f: H \map S^{3}$ so 
that

\begin{enumerate} 

\item $p = p f$ on $N$, i. e. $f$ preserves height and

\item $H \cup (S^{3} - H)$ is a Heegaard splitting of $S^{3}$.

\end{enumerate}
\end{cor}  

\begin{proof} As noted before Proposition \ref{prop:fox}, we may as
well assume that $H$ contains both poles.  The condition on horizontal
disks guarantees, via Proposition \ref{prop:fox}, that the
connectivity graph of $S^{3} - H$ is a tree.  Then Corollary
\ref{cor:reimbed} says there is a height-preserving reimbedding of $H$
so that $S^{3} - H$ is a connected sum of handlebodies.  But since
$\bdd H$ is connected, $S^{3} - H$ is in fact simply a handlebody.
\end{proof}

\section{Knot width} \label{section:width}

For standard definitions about knots in $S^{3}$, see \cite{BZ}, \cite{L} or
\cite{R}.

\begin{defin}
As above, let $p:S^3 \rightarrow R$ be the standard height function
and let $S^{t}$ denote $p^{-1}(t)$, a sphere if $|t| < 1$.  Let $K
\subset S^3$ be a knot in general position with respect to $p$ and
$c_1, \dots, c_n$ be the critical values of $h = p|K$ listed in
increasing order; i.e., so that $c_1 < \dots < c_n$.  Choose $r_1,
\dots, r_{n-1}$ so that $c_i < r_i < c_{i+1}, i = 1, \ldots, n-1$. 
The {\em width of K with respect to $h$}, denoted by $w(K, h)$, is
$\sum_i|K \cap S^{r_{i}}|$.  The {\em width of K}, denoted by $w(K)$,
is the minimum of $w(K', h)$ over all knots $K'$ isotopic to $K$.  We
say that $K$ is in {\em thin position} if $w(K, h) = w(K).$
\end{defin}

We note as an aside that there is an alternative way to calculate
width, inspired by a comment of Clint McCrory.  For the levels 
$r_{i}$ described above, call $r_i$ a {\em
thin level} of $K$ with respect to $h$ if $c_i$ is a maximum value for
$h$ and $c_{i+1}$ is a minimum value for $h$.  Dually $r_i$ is a
{\em thick level} of $K$ with respect to $h$ if $c_i$ is a minimum
value for $h$ and $c_{i+1}$ is a maximum value for $h$.  Since the
lowest critical point of $h$ is a minimum and the highest is a
maximum, there is one more thick level than thin level.

\begin{lemma}
Let $r_{i_1}, \dots, r_{i_k}$ be the thick levels of $K$ and $r_{j_1},
\dots, r_{j_{k-1}}$ the thin levels.  Set $a_{i_l} = |\;K\; \cap
S^{r_{i_l}}|$ and $b_{j_l} = |\;K\; \cap S^{r_{j_l}}|$.  Then \[w(K) = 2
\sum_{l = 1}^k a_{i_l}^2 - 2 \sum_{l = 1}^{k-1} b_{j_l}^2.\]
\end{lemma}

\begin{figure}[tbh] 
\centering
\includegraphics[width=1.5in]{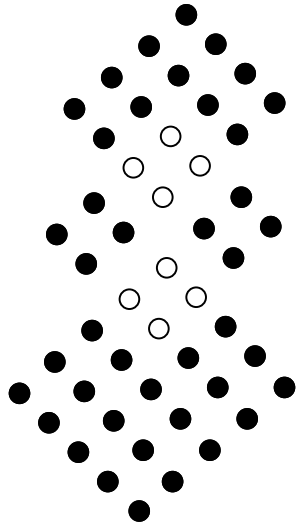}
\caption{Dark dots indicate squares that are added; white dots 
indicate overlap squares that are subtracted} \label{fig:dots}
\end{figure}

\begin{proof}
This can be proven by a direct computation and repeated use of the
Gauss Summation Formula.  It is illustrated in Figure \ref{fig:dots}. 
Each dot represents two points of intersection with a regular level
surface between two critical level surfaces.  For instance, the dots
in Figure \ref{fig:dots} represent the case in which the critical
values, listed from the highest to the lowest are a maximum, maximum,
maximum, maximum, minimum, minimum, maximum, maximum, minimum,
minimum, maximum, maximum, maximum, minimum, minimum, minimum,
minimum, minimum.
\end{proof}

\begin{cor} \label{cor:satellite} Suppose $K$ is a knot in an
unknotted solid torus $W \subset S^{3}$.  Suppose $f: W \map S^{3}$ is
a knotted embedding and $K' = f(K)$.  Then $w(K') \geq w(K)$.
\end{cor}

\begin{proof} Let $p: S^{3} \map R$ be the standard height function. 
Isotope $K'$ so as to minimize its width with respect to this height
function and let $H$ denote the image of $f(W)$ after this isotopy. 
Each generic $2$-sphere $S^{t} = p^{-1}(t)$ intersects $\bdd H$ in a
collection of circles, each of them unknotted since they all lie in
$S^{t}$.  By standard Morse theory, there must be a generic value of
$t$ for which one of the circles $c \subset \bdd H \cap S^{t}$ is
essential in $\bdd H$ and that circle can't be a longitude, since $H$
is a knotted torus.  Hence $c$ must be a meridian circle.  It follows
from Corollary \ref{cor:heegaard} that there is a reimbedding $g$ of
$H$ in $S^{3}$ that preserves height but after which $H$ is unknotted. 
The reimbedding is defined via braid moves on $M = S^{3} - H$; after
perhaps adding a number of Dehn twists to one of the braid moves near
a meridinal boundary component of $P^{t} = M \cap S^{t}$, we can take
this reimbedding to preserve a longitude of $H$.  So in particular,
$g(K')$ is isotopic to $K$ in $S^{3}$ and still has the width of $K'$. 
\end{proof}

Corollary \ref{cor:satellite} can be applied to composite knots, via
the following standard construction.  Let $K = K_1 \# K_2$ be a
composite knot with decomposing sphere $S$.  Then $S^3 - \eta(K \cup
S)$ has two components.  Each of these components is a torus, called a
{\em swallow-follow torus}.  Each of these tori bounds a solid torus in
$S^{3}$ that contains $K$ ; the torus $T_{1}$ whose core is parallel
to $K_{1}$ is said to follow $K_{1}$ and swallow $K_{2}$.  Similarly,
the other torus $T_{2}$ follows $K_{2}$ and swallows $K_{1}$.  The
torus $T_{1}$ exhibits $K$ as a satellite knot of $K_{1}$ with pattern
$K_{2}$, and symmetrically for $T_{2}$.  Therefore, when Corollary
\ref{cor:satellite} is applied to each $T_{i}$ in turn, we get

\begin{cor} \label{cor:sum}
For any two knots $K_1, K_2$, \[w(K_1 \# K_2) \geq max \{ w(K_1),
w(K_2) \} \geq  \frac{1}{2}(w(K_1) + w(K_2)).\]
\end{cor}

Of course the construction can be iterated to give 

\begin{cor}
$w(K_1 \# \dots \# K_n) \geq max\{ w(K_1), \dots, w(K_n) \} \geq
\frac{1}{n}(w(K_1) + \dots + w(K_n)).$
\end{cor}

\begin{proof} For each $K_{i}$ there is a torus that swallows $K_{i}$
and follows the connected sum of the remaining summands.  \end{proof}

It remains to find examples, if any, of knots whose widths degenerate
under connected sum, i.e. knots for which $max \{ w(K_1), w(K_2) \}
\leq w(K_1 \# K_2) < (w(K_1) + w(K_2)) - 2$.


\begin{thebibliography}{30}
    
    
\bibitem[BZ]{BZ} 
G. Burde, H. Zieschang, \underline{Knots}, de Gruyter Studies in
Mathematics 5, Walter de Gruyter \& G., Berlin, 1985, ISBN:
3-11-008675
 
     
\bibitem[Fo]{Fo} R. H.~Fox, \textit{On the imbedding of polyhedra
in $3$-space},  Ann.  of Math. \textbf{49} (1948), 462--470.


\bibitem[G]{G} D. Gabai, \textit{Foliations and the topology of $3$-manifolds. III},
J. Differential Geom. \textbf{26} (1987), 3, 479--536.

\bibitem[L]{L} W.R.B.R. Lickorish, \underline{An introduction to knot theory}, 
Graduate Texts in Mathematics, 175, Springer-Verlag, New York, 1997, ISBN:
0-387-98254-X.

\bibitem[Mo]{Mo} K. Morimoto, \textit{There are knots whose tunnel
numbers go down under connected sum}, Proc.  Am.  Math.  Soc, \textbf{
123} (1995), no.  11, 3527--3532

\bibitem[MS]{MS} K. Morimoto, J. Schultens, \textit{Tunnel numbers of
small knots do not go down under connected sum}, Proc.  Am.  Math. 
Soc, \textbf{ 128} (2000), no.  1, 269--278

\bibitem[RS]{RS} Y.~Rieck, E.~Sedgwick, \textit{Thin position for a
connected sum of small knots}, Algebraic and Geometric Topology
\textbf{2} (2002), 297-309.

\bibitem[R]{R} D. Rolfsen, \underline{Knots and Links} 
Mathematics Lecture Series, No.  7, Publish or Perish,
Inc., Berkeley, Calif., 1976

\bibitem[S]{S} H. Schubert, \textit{\"Uber eine numerische Knoteninvariante},
Math. Z. \textbf{61} (1954), 245-288

\bibitem[Sc1]{Sc1} M. Scharlemann, \textit{Handlebody complements
in the 3-sphere: a remark on a theorm of Fox}, Proc.  Amer.  Math. 
Soc.  \textbf{115} (1992), 1115--1117.

\bibitem[ScSc]{ScSc} M.~Scharlemann, J.~Schultens, {\em Annuli in
generalized Heegaard splittings and degeneration of tunnel number}, 
Math. Ann.  {\bf 317}  (2000),  no. 4, 783--820. 

\bibitem[Sch]{Sc3} J. Schultens, \textit{Additivity of bridge numbers
of knots}, \hspace{2 mm} math.GT/0111032, to appear in Proc. Camb. Phil. 
Soc.

\bibitem[Th]{Th} A. Thompson, personal communication.

\end{thebibliography}
\end{document}